\documentclass[opre,nonblindrev]{informs3modified}

\OneAndAHalfSpacedXI % current default line spacing
\usepackage{endnotes}
\let\footnote=\endnote

\usepackage{natbib}
 \bibpunct[, ]{(}{)}{,}{a}{}{,}%
 \def\bibfont{\small}%
\TheoremsNumberedThrough     % Preferred (Theorem 1, Lemma 1, Theorem 2)
\ECRepeatTheorems

\EquationsNumberedThrough    % Default: (1), (2), ...
\usepackage{bm}
\usepackage{subcaption}
\begin{document}
%%%%%%%%%%%%%%%%

%\setcounter{page}{1}
%\DOI{10.1287/xxxx.0000.0000}%

\RUNAUTHOR{Henry Lam}

\RUNTITLE{Statistical Guarantees via the Empirical DRO}

% \TITLE{Bundling Information Goods of Decreasing Value}
\TITLE{Recovering Best Statistical Guarantees via the Empirical Divergence-based Distributionally Robust Optimization}

\ARTICLEAUTHORS{%
\AUTHOR{Henry Lam}
\AFF{Department of Industrial and Operations Engineering, University of Michigan, Ann Arbor, MI 48109, \EMAIL{khlam@umich.edu}} %, \URL{}}
%\AUTHOR{Marg Arinella}
%\AFF{Institute for Food Adulteration, University of Food Plains, Food Plains, MN 55599, \EMAIL{m.arinella@adult.ufp.edu}}
} % end of the block

\ABSTRACT{%
We investigate the use of distributionally robust optimization (DRO) as a tractable tool to recover the asymptotic statistical guarantees provided by the Central Limit Theorem, for maintaining the feasibility of an expected value constraint under ambiguous probability distributions. We show that using empirically defined Burg-entropy divergence balls to construct the DRO can attain such guarantees. These balls, however, are not reasoned from the standard data-driven DRO framework since by themselves they can have low or even zero probability of covering the true distribution. Rather, their superior statistical performances are endowed by linking the resulting DRO with empirical likelihood and empirical processes. We show that the sizes of these balls can be optimally calibrated using $\chi^2$-process excursion. We conduct numerical experiments to support our theoretical findings.% Enter your abstract
}%

%\KEYWORDS{deterministic inventory theory; infinite linear programming duality;
%  existence of optimal policies; semi-Markov decision process; cyclic schedule}

\KEYWORDS{distributionally robust optimization, empirical likelihood, empirical process, chi-square process, central limit theorem} %\HISTORY{This paper was
%first submitted on April 12, 1922 and has been with the authors for
%83 years for 65 revisions.}

\maketitle
%%%%%%%%%%%%%%%%%%%%%%%%%%%%%%%%%%%%%%%%%%%%%%%%%%%%%%%%%%%%%%%%%%%%%%

% Samples of sectioning (and labeling) in OPRE
% NOTE: (1) \section and \subsection do NOT end with a period
%       (2) \subsubsection and lower need end punctuation
%       (3) capitalization is as shown (title style).
%
%\section{Introduction.}\label{intro} %%1.
%\subsection{Duality and the Classical EOQ Problem.}\label{class-EOQ} %% 1.1.
%\subsection{Outline.}\label{outline1} %% 1.2.
%\subsubsection{Cyclic Schedules for the General Deterministic SMDP.}
%  \label{cyclic-schedules} %% 1.2.1
%\section{Problem Description.}\label{problemdescription} %% 2.

% Text of your paper here

\section{Statistical Motivation of Distributionally Robust Optimization}

We consider an expected value constraint in the form
\begin{equation}
Z_0(x):=E_0[h(x;\xi)]\leq0\label{constraint}
\end{equation}
where $\xi\in\Xi$ is a random object under the probability measure $P_0$, $E_0[\cdot]$ denotes the corresponding expectation, $x\in\Theta\subset\mathbb R^m$ is the decision variable, and $h$ is a known function. The generic constraint \eqref{constraint} has appeared in various applications such as resource allocation (\cite{atlason2004call}), risk management (\cite{palmquist1999portfolio,fabian2008handling}), among others.

In practice, the probability measure $P_0$ is often unknown, but rather is observed via a finite collection of data. Such uncertainty has been considered in the stochastic and the robust optimization literature. Our main goal in this paper is to investigate, in a statistical sense, the \emph{best} data-driven reformulation of \eqref{constraint} in terms of feasibility guarantees.

\subsection{Initial Attempt: Sample Average Approximation}
To define what ``best" means, we start by discussing arguably the most natural attempt for handling \eqref{constraint}, namely the sample average approximation (SAA) (\cite{shapiro2014lectures,wang2008sample,kleywegt2002sample}). Suppose we have i.i.d. data $\xi_1,\ldots,\xi_n$. SAA entails replacing the unknown expectation $Z_0(x)$ with the sample average $(1/n)\sum_{i=1}^nh(x;\xi_i)$, leading to
\begin{equation}
\hat h(x):=\frac{1}{n}\sum_{i=1}^nh(x;\xi_i)\leq0\label{SAA}
\end{equation}
The issue with naively using SAA in this setting is that a solution feasible according to \eqref{SAA} may be mistakeably infeasible for \eqref{constraint}. Since for any $x$ the true mean $Z_0(x)$ can lie above or below its sample average, both with substantial probabilities, the $x$'s close to the boundary of the feasible region according to \eqref{SAA} could, with overwhelming probabilities, be infeasible for the original constraint \eqref{constraint}. %Alternately, letting $\mathcal A=\{x\in\Theta:E_0[h(x;\xi)]\leq0\}$, and $\hat{\mathcal A}=\left\{x\in\Theta:\frac{1}{n}\sum_{i=1}^nh(x;\xi_i)\leq0\right\}$, our discussion suggests equivalently that $P(\mathcal A\subset\hat{\mathcal A})$ can be low.
Consequently, the probability
$$P\left(\hat h(x)\leq0\ \Rightarrow\ Z_0(x)\leq0\right)$$
where $P$ is with respect to the generation of data, can be much lower than an acceptable level.

One way to boost the confidence of SAA is to insert a margin, namely by using the constraint
\begin{equation}
\hat h(x)+\epsilon_n\leq0\label{SAA relaxation}
\end{equation}
This idea has appeared in various contexts (e.g., \cite{wang2008sample,nagaraj2014stochastically}). Choosing $\epsilon_n>0$ suitably can guarantee that
\begin{equation}
P\left(\hat h(x)+\epsilon_n\leq0\ \Rightarrow\ Z_0(x)\leq0\right)\geq1-\alpha\label{conf guarantee}
\end{equation}
where $1-\alpha$ is a prescribed confidence level chosen by the modeler (a typical choice is $\alpha=0.05$). This is achieved by finding $\epsilon_n$ such that
\begin{equation}
P\left(Z_0(x)\leq\hat h(x)+\epsilon_n\text{\ \ for all\ \ }x\in\Theta\right)\geq1-\alpha\label{conf guarantee simplify}
\end{equation}
Such a choice of $\epsilon_n$ can be obtained in terms of the maximal variance of $h(x;\xi)$ over all $x\in\Theta$, and other information such as the diameter of the space $\Theta$ (e.g., \cite{wang2008sample} provides one such choice). %Furthermore, even if these information are available, the margin $\epsilon_n$ can be conservative, in the sense that the LHS of \eqref{conf guarantee simplify} may be far above $1-\alpha$.
%However, As revealed in the literature (e.g., \cite{wang2008sample}),

\subsection{The Statistician's Approach: Confidence Bounds from the Central Limit Theorem}
Though \eqref{conf guarantee simplify} could provide a good feasibility guarantee, the use of one single number $\epsilon_n$ as the margin adjustment may unnecessarily penalize $x$ whose $h(x;\xi)$ bears only a small variation. From a ``classical statistician"'s viewpoint, we adopt a margin adjustment that takes into account the variability of $h(x;\xi)$ at each point of $x$, and at the same time provides a $1-\alpha$ confidence guarantee, by formulating the constraint as
\begin{equation}
\hat h(x)+z\frac{\hat\sigma(x)}{\sqrt n}\leq0\label{CB}
\end{equation}
where $z$ is the critical value of a suitable sampling distribution, and $\hat\sigma(x)$ is an estimate of $\sqrt{Var_0(h(x;\xi))}$ ($Var_0(\cdot)$ denotes the variance under $P_0$), i.e., $\hat\sigma(x)/\sqrt n$ is the standard error. A judicious choice of $z$ can lead to the \emph{asymptotically exact} guarantee
\begin{equation}
\lim_{n\to\infty}P\left(Z_0(x)\leq\hat h(x)+z\frac{\hat\sigma(x)}{\sqrt n}\text{\ \ for all\ \ }x\in\Theta\right)=1-\alpha\label{guarantee exact}
\end{equation}
Without making further assumption on the optimization objective, we set the reformulation \eqref{CB} and the guarantee \eqref{guarantee exact} as our benchmark in this paper, since they stem from the central limit theorem (CLT) widely used in statistics.
%While better (in the sense of being less conservative) reformulation can potentially be constructed by exploiting additional knowledge on the optimization objective and the form of the distribution, in this paper we exclude such information, and view the asymptotically exact guarantee in \eqref{guarantee exact} as the best possible among all constraint reformulations in the form of margin-adjusted SAA, i.e. $(1/n)\sum_{i=1}^nh(x;\xi_i)+\epsilon_n(x)$.

%, we regard the form given by \eqref{CB} that satisfies \eqref{guarantee exact} as the best viewing from the classical asymptotic statistics standpoint. Of course, there is a gap to exactness in any finite sample, and that, given additional structure on the optimization program in hand, one may be able to utilize the information to further improve the feasibility guarantee. But, without  a notion that will be used throughout this paper.

The problem with directly using \eqref{CB} is that (sample) standard deviation is not a tractability-preserving operation, e.g., $\hat\sigma(x)$ may not be convex in $x$ even though the function $h(x;\xi)$ is. Thus the constraint \eqref{CB} can be intractable despite that \eqref{constraint} is tractable. This motivates the investigation of a distributionally robust optimization (DRO) approach, namely, by using
\begin{equation}
\max_{P\in\mathcal U}E_P[h(x;\xi)]\leq0\label{DRO constraint}
\end{equation}
where $E_P[\cdot]$ denotes the expectation under $P$, and $\mathcal U:=\mathcal U(\xi_1,\ldots,\xi_n)$ is an uncertainty set (also known as ambiguity set), calibrated from data, that contains a collection of distributions. As documented in many previous work (e.g., \cite{delage2010distributionally,ben2013robust}), \eqref{DRO constraint} can be made tractable by suitably choosing $\mathcal U$. One central question in this paper is to ask:
\\

\noindent\emph{Is there a tractable choice of $\mathcal U$ that can recover the statistician's asymptotically exact guarantee, namely
\begin{equation}
\lim_{n\to\infty}P\left(Z_0(x)\leq\max_{P\in\mathcal U}E_P[h(x;\xi)]\text{\ \ for all\ \ }x\in\Theta\right)=1-\alpha\label{main guarantee}
\end{equation}
and that
\begin{equation}
\max_{P\in\mathcal U}E_P[h(x;\xi)]\approx\hat h(x)+z\frac{\hat\sigma(x)}{\sqrt n}\ ?\label{mimic CI}
\end{equation}
}

%The rest of this paper is to systematically build a $\mathcal U$ that is very close to answering \eqref{main guarantee}.

\subsection{Data-driven Distributionally Robust Optimization and Statistically ``Good" Uncertainty Sets}\label{sec:good set}
To answer the above question, let us first revisit the common argument in the literature of data-driven DRO. To facilitate discussion, we call an uncertainty set $\mathcal U$ statistically ``good" if it allows
\begin{equation}
\liminf_{n\to\infty}P\left(Z_0(x)\leq\max_{P\in\mathcal U}E_P[h(x;\xi)]\text{\ \ for all\ \ }x\in\Theta\right)\geq1-\alpha\label{good guarantee}
\end{equation}
In contrast, a statistically ``best" uncertainty set in the sense of \eqref{main guarantee} sharpens the inequality in \eqref{good guarantee} to equality.

The data-driven DRO framework provides a general methodology in guaranteeing \eqref{good guarantee}. First, one calibrates an uncertainty set $\mathcal U$ from data so that it contains the true distribution with probability $1-\alpha$, namely $P(P_0\in\mathcal U)\geq1-\alpha$. Note that since $P_0\in\mathcal U$ implies that $Z_0(x)=E_0[h(x;\xi)]\leq\max_{P\in\mathcal U}E_P[h(x;\xi)]$ for all $x$, we have
\begin{equation}
P\left(Z_0(x)\leq\max_{P\in\mathcal U}E_P[h(x;\xi)]\text{\ \ for all\ \ }x\in\Theta\right)\geq P(P_0\in\mathcal U)\geq1-\alpha\label{good guarantee argument}
\end{equation}
Similarly, a set $\mathcal U$ constructed with the asymptotic property $\liminf_{n\to\infty}P(P_0\in\mathcal U)\geq1-\alpha$ guarantees that \eqref{good guarantee} holds and, in fact, so is the stronger guarantee
$$\liminf_{n\to\infty}P\left(\min_{P\in\mathcal U}E_P[h(x;\xi)]\leq Z_0(x)\leq\max_{P\in\mathcal U}E_P[h(x;\xi)]\text{\ \ for all\ \ }x\in\Theta\right)\geq1-\alpha$$

Thus, good uncertainty sets can be readily created as confidence regions for $P_0$. Constructing these confidence regions and their tractability have been substantially investigated. A non-exhaustive list includes moment and deviation-type constraints (\cite{delage2010distributionally,goh2010distributionally,wiesemann2014distributionally}), Wasserstein balls (\cite{esfahani2015data,gao2016distributionally}), $\phi$-divergence balls (\cite{ben2013robust}), likelihood-based (\cite{Wang2015}) and goodness-of-fit-based regions (\cite{bertsimas2014robust}). Recently, \cite{gupta2015near} further investigates the smallest of such confidence regions as a baseline to measure the degree of conservativeness of a given uncertainty set. %For convenience, we call the use of

\subsection{Our Contributions}
Despite the availability of all the good uncertainty sets, finding the statistically best one in the sense of \eqref{main guarantee} has not been addressed in the literature. In this paper, we construct an uncertainty set that is close to the best (the meaning of ``close to" will be apparent in our later exposition) by leveraging one of the good sets, namely the Burg-entropy divergence ball.

Intriguingly, the way we construct these balls, and the associated statistical explanation, is completely orthogonal to the standard data-driven DRO framework discussed above. These balls are empirically defined (as we will explain in detail) and do not have any interpretation as confidence regions by themselves. In fact, they have low, or even zero, probability of covering the true distribution. Yet the resulting DRO has the best statistical performances among all DRO formulations. This disentanglement between set coverage and ultimate performance can be explained by a duality relation between our resulting DRO and the empirical likelihood theory, a connection that has been briefly discussed in a few previous work (e.g., \cite{Wang2015,lam2015quantifying}) but not been fully exploited as far as we know.
%reveals the fact that merely choosing a good confidence region as an uncertainty set does not automatically translate to good statistical properties; more importantly, the converse holds that an invalid confidence region turns out to give the best ultimate guarantees.

Importantly, through setting up such a connection, we study optimal calibration of the sizes of these sets by using a generalization of $\chi^2$-quantiles that involves the excursion of so-called $\chi^2$-processes. As a by-product, our proposed method also resolves some technical challenges reported in the previous literature in calibrating divergence balls (e.g., \cite{jiang2012data,esfahani2015data}). More precisely, since divergence is only properly defined between absolutely continuous distributions, it has been suggested, in the case of continuous distributions, that one needs to construct the ball using kernel estimation of density and the divergence, which is statistically challenging, or resorting to a parametric framework. The approach we take here, on the other hand, bypasses these issues.

To summarize, our main contributions of this paper are:\emph{
\begin{enumerate}
\item We systematically build an uncertainty set that, in a precise sense, is close to recovering the guarantees \eqref{main guarantee} and \eqref{mimic CI} provided by the CLT.
\item In doing so, we expand the view on the meaning of uncertainty sets beyond the notion of confidence regions, by showing that our empirical Burg-entropy divergence ball recovers the best guarantees despite being a low or zero-coverage set. This is achieved through connecting the dual of the resulting DRO with the empirical likelihood theory.
\item To achieve our claimed guarantees, we study an approach to optimally calibrate the sizes of these balls using quantiles of $\chi^2$-process excursion.
\item As a by-product, our approach resolves the technical difficulties in enforcing absolute continuity when calibrating divergence balls that are raised in previous works in data-driven DRO.
\end{enumerate}
}

%We note that the viewpoint taken in this paper is purely statistical. Modeling a given problem may involve other important dimensions of considerations such as risks and robustness (i.e., optimality and feasibility performances of the solution under perturbation of parameters) etc.. The latter is part of the original motivation for using robust optimization (see, e.g., \cite{ben2009robust} and \cite{bertsimas2011theory}). These considerations are out of the scope of this work.
Finally, while the viewpoint taken by this paper is primarily statistical, we mention that there are other valuable perspectives in the DRO literature motivated from risk or tractability considerations (see, e.g., the survey \cite{gabrel2014recent}); these are, however, beyond the scope of this work.

The rest of this paper is organized as follows. Section \ref{sec:empircal uncertainty set} motivates our proposed uncertainty sets. Section \ref{sec:theory} presents methods to calibrate their sizes and the theoretical explanation of their statistical performances. Section \ref{sec:numerics} shows results of our numerical experiments. Section \ref{sec:conlusion} concludes and discusses future directions. Section \ref{sec:proofs} provides all the proofs. Appendices \ref{sec:EP} and \ref{sec:thms} list some auxiliary concepts and theorems.

%These contributions center around what we call \emph{empirical} divergence-based uncertainty sets, which we discuss next.
%By ``confidence guarantee" we mean the probability of the decision, a.k.a. solution, being feasible is at least $1-\alpha$, where $\alpha$ is a prescribed level set by the modeler (a typical choice is $\alpha=0.05$). By ``best" here we mean, without further knowledge on the structure of the umbrella optimization, having the least level of conservativeness. This latter terminology is more subtle, and we shall discuss further momentarily.
%
%In order to clarify and explain our goal further, let us begin with the most natural attempt
%

\section{Towards the Empirical DRO}\label{sec:empircal uncertainty set}
We first review some background in divergence-based inference and how to use it to create confidence regions for probability distributions in Section \ref{sec:DI}. Through a preliminary numerical investigation in Section \ref{sec:initial}, we motivate and present, in Section \ref{sec:empirical divergence ball}, the \emph{empirical divergence ball} as our main tool.

\subsection{Divergence-based Inference and Confidence Regions}\label{sec:DI}
A $\phi$-divergence ball is in the form
\begin{equation}
\mathcal U=\{P\in\mathcal P_Q:D_\phi(P,Q)\leq\eta\}\label{divergence}
\end{equation}
where
$$D_\phi(P,Q)=\int\phi\left(\frac{dP}{dQ}\right)dQ$$
for some baseline distribution $Q$ and suitable function $\phi(\cdot)$, and $dP/dQ$ is the likelihood ratio given by the Radon-Nikodym derivative between $P$ and $Q$. The latter is well-defined only for $P$ within $\mathcal P_Q$, the set of all distributions absolutely continuous with respect to $Q$. The function $\phi:\mathbb R^+\to\mathbb R$ is convex and satisfies $\phi(1)=0$.

Suppose the random variable $\xi$ lies on a finite discrete support $\{s_1,\ldots,s_k\}$. One way to construct a statistically good divergence ball is as follows (\cite{ben2013robust}). Set the baseline distribution as the histogram of the i.i.d. data given by $\hat{\mathbf p}=(\hat p_i)_{i=1,\ldots,k}$, where $\hat p_i=n_i/n$, $n_i$ is the counts on support $s_i$, and $n$ is the total sample size. The divergence ball \eqref{divergence} can be written as
\begin{align}
\mathcal U&=\{\mathbf p\in\mathcal P_{\hat{\mathbf p}}:D_\phi(\mathbf p,\hat{\mathbf p})\leq\eta\}\notag\\
&=\left\{(p_1,\ldots,p_k):\sum_{i=1}^k\hat p_i\phi\left(\frac{p_i}{\hat p_i}\right)\leq\eta,\ \sum_{i=1}^kp_i=1,\ p_i\geq0\text{\ for all\ }i=1,\ldots,k\right\}\label{divergence discrete}
\end{align}
Under twice continuous differentiability condition on $\phi$, the theory of divergence-based inference (\cite{pardo2005statistical}) stipulates that
$$\frac{2n}{\phi''(1)}D_\phi(\mathbf p,\hat{\mathbf p})\Rightarrow\chi^2_{k-1}\text{\ \ as $n\to\infty$}$$
where $\chi^2_{k-1}$ is the $\chi^2$-distribution with degree of freedom $k-1$, and ``$\Rightarrow$" denotes convergence in distribution. This implies that taking $\eta=\frac{\phi''(1)}{2n}\chi^2_{k-1,1-\alpha}$ in \eqref{divergence discrete}, where $\chi^2_{k-1,1-\alpha}$ is the $1-\alpha$ quantile of $\chi^2_{k-1}$, forms an uncertainty set $\mathcal U$ that contains the true distribution with probability asymptotically $1-\alpha$. This in turn implies that $\mathcal U$ is a good uncertainty set satisfying \eqref{good guarantee}.

For instance, $\phi(x)=(x-1)^2$ yields the $\chi^2$-distance, and setting $\eta$ at $\chi^2_{k-1,1-\alpha}/n$ results in the confidence region associated with the standard $\chi^2$ goodness-of-fit test for categorical data (\cite{agresti2011categorical}). On the other hand, $\phi(x)=-\log x+x-1$ yields the Burg-entropy (or the Kullback-Leibler) divergence (\cite{kullback1951information}), and $\eta$ in this case should be set at $\chi^2_{k-1,1-\alpha}/(2n)$. Since the Burg-entropy divergence is important in our subsequent discussion, for convenience, we denote its divergence ball as
\begin{equation}
\mathcal U_{Burg}=\left\{(p_1,\ldots,p_k):-\sum_{i=1}^k\hat p_i\log\frac{p_i}{\hat p_i}\leq\frac{\chi^2_{k-1,1-\alpha}}{2n},\ \sum_{i=1}^kp_i=1,\ p_i\geq0\text{\ for all\ }i=1,\ldots,k\right\}\label{Burg}
\end{equation}
From the discussion above, $\mathcal U_{Burg}$ is a good uncertainty set and moreover satisfies
\begin{equation}
\lim_{n\to\infty}P\left(P_0\in\mathcal U_{Burg}\right)=1-\alpha\label{Burg guarantee}
\end{equation}
for a finite discrete true distribution $P_0$.

The computational tractability of divergence balls has been studied in depth in \cite{ben2013robust}, who reformulate $\max_{P\in\mathcal U}E_P[h(x;\xi)]$ in terms of the conjugate function of $\phi$ and propose efficient optimization algorithms. Because of this we will not drill further on tractability and instead refer interested readers therein.

%However, there remains one crucial gap in applying these results in data-driven environments: While calibration methods for $\eta$ in \eqref{divergence} has been established when the random object $\xi$ is discrete, generalizing these methods to continuous $\xi$ is difficult as the goodness-of-fit statistic underpinning the former, which utilizes the empirical distribution as the baseline, does not hold for continuous data since the truths in such cases are not absolutely continuous to any empirical distributions (as pointed out by several authors, e.g. REF). The only known available resort involves estimation of divergences, which is statistically challenging as they consist of functionals of densities (REF).

\subsection{An Initial Numerical Investigation on Coverage Accuracy}\label{sec:initial}
To get a sense of the coverage performance provided by $\mathcal U_{Burg}$, we run an experiment on estimating $Z_0(x)=E_0[h(x;\xi)]$, where we set $h$ as
\begin{equation}
h(x;\xi)=-v\min(x,\xi)-s(x-\xi)^++l(\xi-x)^++cx+\rho\label{h example}
\end{equation}
with $v=10$, $s=5$, $l=4$, $c=3$, and $\rho=40$. This function $h$ is adapted from the example in Section 6.3 in \cite{ben2013robust}. As an application, \eqref{h example} can represent the loss amount in excess of the threshold $\rho$ for a newsvendor. In this case, $v$ is the selling price per unit, $s$ the salvage value per unit, $l$ the shortage cost per unit, $c$ the cost per unit, $\xi$ a random demand, and $x$ the quantity to order.

For now, let us fix the solution at $x=30$ (so it is purely about estimating $Z_0(30)$). We set the random variable $\xi$ as an exponential random variable with mean $20$ that is discretized uniformly over a $k$-grid on the interval $[0,50]$, or more precisely,
\begin{equation}
\begin{array}{l}
P\left(\xi=\frac{50j}{k}\right)=P\left(\frac{50(j-1)}{k}<Exp\left(\frac{1}{20}\right)<\frac{50j}{k}\right)\text{\ \ for\ }j=1,\ldots,k-1\\
P\left(\xi=50\right)=P\left(Exp\left(\frac{1}{20}\right)>\frac{50(k-1)}{k}\right)\end{array}\label{exp specification}
\end{equation}
We repeat $1,000$ times:
\begin{enumerate}
\item Simulate $n$ i.i.d. data $\xi_1,\ldots,\xi_n$ from the $k$-discretized $Exp(1/20)$.
\item Construct $\mathcal U_{Burg}$, and compute $\min_{\mathbf p\in\mathcal U_{Burg}}E_{\mathbf p}[h(x;\xi)]$ and $\max_{\mathbf p\in\mathcal U_{Burg}}E_{\mathbf p}[h(x;\xi)]$ % where $\mathcal U_{DI}$ is the Burg-entropy divergence ball
%$$\mathcal U=\left\{(p_1,\ldots,p_k):-\sum_{i=1}^k\hat p_i\log\frac{p_i}{\hat p_i}\leq\frac{\chi^2_{k-1,1-\alpha}}{2n},\ \sum_{i=1}^kp_i=1,\ p_i\geq0\text{\ for all\ }i=1,\ldots,k\right\}$$
with $\alpha=0.05$.
\item Output $I\left(\min_{\mathbf p\in\mathcal U_{Burg}}E_{\mathbf p}[h(x;\xi)]\leq Z_0(x)\leq\max_{\mathbf p\in\mathcal U_{Burg}}E_{\mathbf p}[h(x;\xi)]\right)$, where $Z_0(x)$ is the true quantity calculable in closed-form, and $I(\cdot)$ is the indicator function.
\end{enumerate}
We then output the point estimate and the $95\%$ confidence interval (CI) of the coverage probability from the $1,000$ replications.

Step 2 above is carried out by using duality and numerically solving
\begin{align*}
\min_{\mathbf p\in\mathcal U_{Burg}}E_{\mathbf p}[h(x;\xi)]&=\max_{\lambda\geq0,\gamma}\sum_{i=1}^n\frac{\lambda}{n}\log\left(1-\frac{-h(\xi_i)+\gamma}{\lambda}\right)-\lambda\eta+\gamma\\
\max_{\mathbf p\in\mathcal U_{Burg}}E_{\mathbf p}[h(x;\xi)]&=\min_{\lambda\geq0,\gamma}-\sum_{i=1}^n\frac{\lambda}{n}\log\left(1-\frac{h(\xi_i)+\gamma}{\lambda}\right)+\lambda\eta-\gamma
\end{align*}
where $-0\log(1-t/0):=0$ for $t\leq0$ and $-0\log(1-t/0):=\infty$ for $t>0$ (see \cite{ben2013robust}).

Table \ref{table:DRO} shows the estimates of coverage probabilities for different support size $k$. The sample size for $\xi$ is $n=30$. The coverage probabilities are all greater than $95\%$, showing correct statistical guarantees. However, more noticeable is that they are all higher than $99\%$, and are consistently close to $100\%$ for $k=10$ or above, thus leading to severe over-coverage. Note that this phenomenon occurs despite that $\mathcal U_{Burg}$ has asymptotically exactly $1-\alpha$ probability of covering the true distribution as guaranteed in \eqref{Burg guarantee}.

\begin{table}[!htb]
  \begin{subtable}{0.33\textwidth}
   \centering
    {\footnotesize\begin{tabular}{c|c|c}
    $k$&Cover.&$95\%$ C.I. of\\
    &Prob.&Cover. Prob.\\
    \hline
    5&99.6\%&(99.3\%, 99.9\%)\\
10&100.0\%&(100.0\%, 100.0\%)\\
15&100.0\%&(100.0\%, 100.0\%)\\
20&100.0\%&(100.0\%, 100.0\%)
    \end{tabular}}
      \subcaption{DRO with Burg-ball of size $\chi^2_{k-1,0.95}/(2n)$}
    %     \captionof{table}{Coverage probabilities for distributions of different support sizes using data-driven DRO}
      \label{table:DRO}
  \end{subtable}
  %\hspace{.3cm}
  \begin{subtable}{0.33\textwidth}
    \centering
    {\footnotesize\begin{tabular}{c|c|c}
    $k$&Cover.&$95\%$ C.I. of\\
    &Prob.&Cover. Prob.\\
    \hline
        5&94.5\%&(93.3\%, 95.7\%)\\
10&95.0\%&(94.0\%, 96.2\%)\\
15&95.3\%&(94.2\%, 96.4\%)\\
20&94.8\%&(93.7\%, 96.0\%)
    \end{tabular}}
    \subcaption{Standard CLT}
     %     \captionof{table}{Coverage probabilities for distributions of different support sizes using standard CLT}
      \label{table:CLT}
    \end{subtable}
  \begin{subtable}{0.33\textwidth}
     \centering
    {\footnotesize\begin{tabular}{c|c|c}
    $k$&Cover.&$95\%$ C.I. of\\
    &Prob.&Cover. Prob.\\
    \hline
        5&94.1\%&(92.9\%, 95.3\%)\\
10&94.4\%&(93.2\%, 95.6\%)\\
15&95.4\%&(94.3\%, 96.5\%)\\
20&95.3\%&(94.2\%, 96.4\%)
    \end{tabular}}
     \subcaption{DRO with Burg-ball of size $\chi^2_{1,0.95}/(2n)$}
    %     \captionof{table}{Coverage probabilities for distributions of different support sizes using empirical DRO}
      \label{table:empirical}
       \end{subtable}
       \caption{Coverage probabilities for different methods and support sizes for discrete distributions}
  \end{table}

As a comparison, we repeat the experiment, but this time checking the coverage of the standard $95\%$ CI generated from the CLT
$$\left[\hat h(30)-z_{1-\alpha/2}\frac{\hat\sigma(30)}{\sqrt n},\hat h(30)+z_{1-\alpha/2}\frac{\hat\sigma(30)}{\sqrt n}\right]$$
where $\hat h(x)=\frac{1}{n}\sum_{i=1}^nh(x;\xi_i)$, $\hat\sigma^2(x)=\frac{1}{n-1}\sum_{i=1}^n(h(x;\xi_i)-\bar h)^2$, and $z_{1-\alpha/2}$ is the $(1-\alpha/2)$-quantile of standard normal distribution. Table \ref{table:CLT} shows that, unlike data-driven DRO, the coverage probabilities are now very close to $95\%$, regardless of the values of $k$. This result is, of course, as predicted by the CLT.

To investigate the source of inferiority in the data-driven DRO approach, we shall interpret the degree of freedom in the $\chi^2$-distribution from another angle: In maximum likelihood theory (\cite{cox1979theoretical}), the degree of freedom in the limiting $\chi^2$-distribution of the so-called log-likelihood ratio is equal to the number of effective parameters to be estimated. In our experiment, this number is \emph{one}, because we are only interested in estimating a single quantity $Z_0(30)$. Indeed, Table \ref{table:empirical} shows that the coverage probabilities of DRO, using the quantile of $\chi^2_1$ instead of $\chi^2_{k-1}$, are equally competitive as the CLT approach. This motivates us to propose our key definition of uncertainty set next.

\subsection{The Empirical Divergence Ball}\label{sec:empirical divergence ball}
Given i.i.d. data $\xi_1,\ldots,\xi_n$, we define the \emph{empirical Burg-entropy divergence ball} as
\begin{equation}
\mathcal U_n(\eta)=\left\{\mathbf w=(w_1,\ldots,w_n):-\frac{1}{n}\sum_{i=1}^n\log(nw_i)\leq\eta,\ \sum_{i=1}^nw_i=1,\ w_i\geq0\text{\ for all\ }i=1,\ldots,n\right\}\label{empirical uncertainty}
\end{equation}
where $\mathbf w$ is the probability weight vector on the $n$ support points from data (some possibly with the same values). The set \eqref{empirical uncertainty} is well-defined whether the distribution of $\xi$ is discrete or continuous. It is a Burg-entropy divergence ball centered at the empirical distribution, with radius $\eta>0$, pretending that the support of the distribution is solely on the data. For convenience, we call the corresponding DRO over the empirical divergence ball as the \emph{empirical DRO}.

%the Burg-entropy divergence (REF) between $\mathbf w$ and the empirical distribution from data. We choose to focus on Burg-entropy divergence since, as we will see, it has a natural interpretation in the theory of the EL method that we will use heavily. %We call \eqref{empirical uncertainty} the \emph{empirical divergence-based uncertainty set}, and \eqref{empirical opt} the \emph{empirical

The discussion in Section \ref{sec:initial} suggests to put $\eta=\chi^2_{1,1-\alpha}/(2n)$. One intriguing observation is that $\mathcal U_n(\chi^2_{1,1-\alpha}/(2n))$ under-covers the true probability distribution. This can be seen by noting that, in the discrete case, $\mathcal U_n(\chi^2_{1,1-\alpha}/(2n))$ is equivalent to $\mathcal U_{Burg}$ except that $\chi^2_{k-1,1-\alpha}$ in its definition is replaced by $\chi^2_{1,1-\alpha}$ (see Proposition \ref{EL discrete} later). Since $\mathcal U_{Burg}$ is asymptotically exact in providing $1-\alpha$ coverage for the true distribution, and that $\chi^2_{1,1-\alpha}<\chi^2_{k-1,1-\alpha}$, $\mathcal U_n(\chi^2_{1,1-\alpha}/(2n))$ must be asymptotically under-covering. What is more, in the continuous case, the empirical distribution is singular with respect to the true distribution. Thus $\mathcal U_n(\chi^2_{1,1-\alpha}/(2n))$ has as low as \emph{zero} coverage. Clearly, the performance of the empirical uncertainty set cannot be reasoned using the standard data-driven DRO framework discussed in Section \ref{sec:good set}.

We close this section with Table \ref{table:cont}, which shows additional experimental results for the same example as above, this time the data being generated from the \emph{continuous} distribution $Exp(1/20)$. As we can see in Table \ref{table:empirical cont}, the coverages using $\mathcal U_n(\chi^2_{1,1-\alpha}/(2n))$ are maintained at close to $95\%$ when $n=40$ or above. As a comparison, Table \ref{table:CLT cont} shows that the standard CLT performs similarly as the empirical DRO (except that it tends to over-cover instead of under-cover when $n$ is small). Note that, unlike the discrete case, there is no well-defined choice of $k$ in this setting.

\begin{table}[!htb]
  \begin{subtable}{0.5\textwidth}
      \centering
    \begin{tabular}{c|c|c}
    $n$&Cover.&$95\%$ C.I. of\\
    &Prob.&Cover. Prob.\\
    \hline
    20&91.9\%&(90.5\%, 93.3\%)\\
30&92.8\%&(91.5\%, 94.1\%)\\
%40&94.3\%&(93.1\%, 95.5\%)\\
50&94.5\%&(93.3\%, 95.7\%)\\
%60&94.5\%&(93.3\%, 95.7\%)\\
80&94.4\%&(93.2\%, 95.6\%)
    \end{tabular}
    \subcaption{Empirical DRO with ball size $\chi^2_{1,0.95}/(2n)$}
    %     \captionof{table}{Coverage probabilities for distributions of different support sizes using empirical DRO}
      \label{table:empirical cont}
       \end{subtable}
       \begin{subtable}{0.5\textwidth}
    \centering
    \begin{tabular}{c|c|c}
    $n$&Cover.&$95\%$ C.I. of\\
    &Prob.&Cover. Prob.\\
    \hline
     20&96.1\%&(95.1\%, 97.1\%)\\
30&96.4\%&(95.4\%, 97.4\%)\\
%40&94.7\%&(93.5\%, 95.9\%)\\
50&94.3\%&(93.1\%, 95.5\%)\\
%60&94.9\%&(93.8\%, 96.0\%)\\
80&96.4\%&(95.4\%, 97.4\%)
    \end{tabular}
     \subcaption{Standard CLT}
    %     \captionof{table}{Coverage probabilities for distributions of different support sizes using standard CLT}
      \label{table:CLT cont}
    \end{subtable}
      \caption{Coverage probabilities for different methods and sample sizes for continuous distributions}
 \label{table:cont}
  \end{table}

%The next section presents the theory to explain the performance of the empirical divergence ball above, and from there we will develop the calibration method in achieving the best guarantees in \eqref{main guarantee} and \eqref{mimic CI}.

\section{Statistical Guarantees}\label{sec:theory}
We present our theoretical justification in two subsections. Section \ref{sec:EL} first connects the dual of the empirical DRO with the empirical likelihood (EL) method. Sections \ref{sec:excursion}, \ref{sec:PNLRP} and \ref{sec:Euler} elaborate this connection to develop the calibration method for the radius $\eta$ in the empirical divergence ball, via estimating the excursion of $\chi^2$-processes. We defer all proofs to Appendix \ref{sec:proofs}. %The first key observation is that the dual of the empirical DRO is essentially , which we will discuss next.

Throughout our exposition, ``$\Rightarrow$" denotes weak convergence (or convergence in distribution), ``a.s." abbreviates ``almost surely", and ``ev." abbreviates ``eventually".

\subsection{The Empirical Likelihood Method}\label{sec:EL}
The EL method, first proposed by Owens (\cite{owen1988empirical,owen2001empirical}), can be viewed as a nonparametric counterpart of maximum likelihood theory. Given a set of i.i.d. data $\xi_1,\ldots,\xi_n$, one can view the empirical distribution, formed by putting probability weight $1/n$ on each data point, as a nonparametric maximum likelihood in the following sense. We define the nonparametric likelihood of any distributions supported on the data as
\begin{equation}
\prod_{i=1}^nw_i\label{nonparametric likelihood}
\end{equation}
where $\mathbf w=(w_1,\ldots,w_n)\in\mathcal P_n$ is any probability vector on $\{\xi_1,\ldots,\xi_n\}$. Then the likelihood of the empirical distribution, given by
\begin{equation}
\prod_{i=1}^n\frac{1}{n}\label{nonparametric MLE}
\end{equation}
maximizes \eqref{nonparametric likelihood}. This observation can be easily verified by a simple convexity argument. Moreover, \eqref{nonparametric MLE} still maximizes even if one considers other distributions that are not only supported on the data, since these distributions would have $\sum_{i=1}^nw_i<1$, making \eqref{nonparametric likelihood} even smaller.

The key of EL is a nonparametric analog of the celebrated Wilks' Theorem (\cite{cox1979theoretical}), stating the convergence of the so-called logarithmic likelihood ratio to $\chi^2$-distribution. In the EL framework, the nonparametric likelihood ratio is defined as the ratio between any nonparametric likelihood and the maximum likelihood, given by
$$\prod_{i=1}^n\frac{w_i}{1/n}=\prod_{i=1}^n(nw_i)$$
%To establish a nonparametric analog of Wilks' Theorem, one needs to appropriately define a ``parameter", and consequently a ``degree of freedom". Taking a ``goal-driven" perspective, suppose now that we have some performance measure $\mu=E[h(\xi)]$ that we would like to estimate. We treat $\mu$ as a parameter, and consider the so-called nonparametric profile likelihood ratio
To carry out inference we need to specify a quantity of interest to be estimated. Suppose we are interested in estimating $\mu_0=E_0[g(\xi)]$ for some function $g(\cdot)$, where $E_0[\cdot]$ is the expectation with respect to the true distribution generating the data (and similarly, $Var_0(\cdot)$ denotes its variance). The EL method utilizes the \emph{profile nonparametric likelihood ratio}
\begin{equation}
R(\mu)=\max\left\{\prod_{i=1}^n(nw_i):\sum_{i=1}^ng(\xi_i)w_i=\mu,\ \sum_{i=1}^nw_i=1,\ w_i\geq0\text{\ for all\ }i=1,\ldots,n\right\}\label{profile LR}
\end{equation}
where the likelihood ratios are ``profiled" according to the value of $\sum_{i=1}^ng(\xi_i)w_i$. With this definition, we have: %(this technique is used commonly in nested hypothesis tests in the parametric regime; REF). With these developments, the main result is the empirical likelihood theorem (REF):

\begin{theorem}[The Empirical Likelihood Theorem; \cite{owen1988empirical}]
Let $\xi_1,\ldots,\xi_n\in\Xi$ be i.i.d. data under $P_0$. Let $\mu_0=E_0[g(\xi)]<\infty$, and assume that $0<Var_0(g(\xi))<\infty$. Then
\begin{equation}
-2\log R(\mu_0)\Rightarrow\chi^2_1\text{\ \ as $n\to\infty$}\label{Wilks EL}
\end{equation}
where $-2\log R(\mu_0)$ is defined as $\infty$ if there is no feasible solution in defining $R(\mu_0)$ in \eqref{profile LR}.\label{ELT}
\end{theorem}
The degree of freedom 1 in the limiting $\chi^2$-distribution in \eqref{Wilks EL} counts the number of effective parameters, which is only $\mu_0$ in this case. %This is a consequence of a univariate CLT that is embedded in the asymptotic theory of $-2\log R(\mu_0)$.

Phrasing in terms of our problem setup, we define
\begin{equation}
R(x;Z)=\max\left\{\prod_{i=1}^n(nw_i):\sum_{i=1}^nh(x;\xi_i)w_i=Z(x),\ \sum_{i=1}^nw_i=1,\ w_i\geq0\text{\ for\ }i=1,\ldots,n\right\}\label{EL process1}
\end{equation}
and hence
\begin{eqnarray*}
&&-2\log R(x;Z)\\
&=&\min\left\{-2\sum_{i=1}^n\log(nw_i):\sum_{i=1}^nh(x;\xi_i)w_i=Z(x),\ \sum_{i=1}^nw_i=1,\ w_i\geq0\text{\ for all\ }i=1,\ldots,n\right\}
\end{eqnarray*}
From Theorem \ref{ELT}, we conclude $P(-2\log R(x;Z_0)\leq\chi^2_{1,1-\alpha})\to1-\alpha$ as $n\to\infty$ for a fixed $x$. The important implication of Theorem \ref{ELT} arises from a duality relation between $-2\log R(x;Z_0)$ and the optimal values of the empirical DRO, in the sense that $-2\log R(x;Z_0)\leq\kappa$
if and only if
$$\min_{\mathbf w\in\mathcal U_n(\kappa/(2n))}\sum_{i=1}^nh(x;\xi_i)w_i\leq Z_0(x)\leq\max_{\mathbf w\in\mathcal U_n(\kappa/(2n))}\sum_{i=1}^nh(x;\xi_i)w_i$$
where $\mathcal U_n(\eta)$ is the empirical divergence ball defined in \eqref{empirical uncertainty}.
This implies:
\begin{theorem}
Fix $x\in\Theta$, and let $\xi_1,\ldots,\xi_n\in\Xi$ be i.i.d. data under $P_0$. Assume that $0<Var_0(h(x;\xi))<\infty$, and $Z_0(x)=E_0[h(x;\xi)]<\infty$. We have
\begin{equation}
\lim_{n\to\infty}P\left(\underline Z_n(x)\leq Z_0(x)\leq\overline Z_n(x)\right)=1-\alpha\label{exact asymptotic}
\end{equation}
where
\begin{align}
\underline Z_n(x)&=\min_{\mathbf w\in\mathcal U_n\left(\chi^2_{1,1-\alpha}/(2n)\right)}\sum_{i=1}^nh(x;\xi_i)w_i\label{min single}\\
\overline Z_n(x)&=\max_{\mathbf w\in\mathcal U_n\left(\chi^2_{1,1-\alpha}/(2n)\right)}\sum_{i=1}^nh(x;\xi_i)w_i\label{max single}
\end{align}
% with
%\begin{equation}
%\eta=\frac{\chi^2_{1,1-\alpha}}{2n}\label{basic calibration}
%\end{equation}
%where $\chi^2_{1,1-\alpha}$ is the $(1-\alpha)$-quantile of $\chi^2_1$, the $\chi^2$ distribution with degree of freedom 1.
\label{EL basic}
\end{theorem}

Next, we argue that, in the discrete case, the empirical DRO given by $\underline Z_n(x)$ and $\overline Z_n(x)$ reduces to the standard divergence-based DRO given by $\max/\min_{\mathbf p\in\mathcal U_{Burg}}E_{\mathbf p}[h(x;\xi)]$, except that the degree of freedom in the $\chi^2$-quantile is replaced by 1. This explains the experimental results in Section \ref{sec:initial}.
\begin{proposition}
Fix $x\in\Theta$. When $\xi$ is discrete on the support set $\{s_1,\ldots,s_k\}$, $\underline Z_n(x)$ and $\overline Z_n(x)$ defined in \eqref{min single} and \eqref{max single} are equal to $\min_{\mathbf p\in\mathcal U_{Burg}'}E_{\mathbf p}[h(x;\xi)]$ and $\max_{\mathbf p\in\mathcal U_{Burg}'}E_{\mathbf p}[h(x;\xi)]$ respectively, where
\begin{equation}
\mathcal U_{Burg}'=\left\{(p_1,\ldots,p_k):-\sum_{i=1}^k\hat p_i\log\frac{p_i}{\hat p_i}\leq\frac{\chi^2_{1,1-\alpha}}{2n},\ \sum_{i=1}^kp_i=1,\ p_i\geq0\text{\ for all\ }i=1,\ldots,k\right\}\label{Burg EL}
\end{equation}
and $\hat p_i=n_i/n$, the proportion of data falling onto $s_i$.\label{EL discrete}
\end{proposition}

We complement Theorem \ref{EL basic} with a consistency result:
\begin{theorem}
Under the same conditions in Theorem \ref{EL basic}, for any fixed $x\in\Theta$, both $\underline Z_n(x)\stackrel{a.s.}{\to}Z_0(x)$ and $\overline Z_n(x)\stackrel{a.s.}{\to}Z_0(x)$ as $n\to\infty$.\label{consistency simple}
\end{theorem}

We note that, in the data-driven DRO framework,  if $\xi$ is continuous, the absolute continuity condition requires a divergence ball to center at a continuous distribution to have any chance of containing the true distribution. This observation has been pointed out by several authors (e.g., \cite{jiang2012data,esfahani2015data}) and forces the use of kernel density estimators to set the baseline. Unless one assumes a parametric framework, calibrating the ball radius requires nonparametric divergence estimation, which involves challenging statistical analyses on bandwidth tuning and loss of estimation efficiency (e.g., \cite{moon2014multivariate,nguyen2007estimating,pal2010estimation}). The empirical DRO based on the EL framework cleanly bypasses these issues.

Our discussion in this subsection is also related to likelihood robust optimization studied in \cite{Wang2015}, which also discusses EL as well as other connections such as Bayesian statistics. \cite{Wang2015} focuses on finite discrete distributions. The work \cite{lam2015quantifying} also investigates EL, among other techniques like the bootstrap, in constructing confidence bounds for the optimal values of stochastic programs. However, none of these formalizes the connection, or more precisely, the \emph{disconnection} between set coverage and the statistical performance of DRO. As our next subsection shows, this formalization is important in capturing a statistical price to attain our best guarantee in \eqref{main guarantee}. This will be our focus next.

%We will prove a generalized version of this theorem later.

%Despite the apparent simplicity in calibrating the empirical uncertainty set in an estimation problem, namely when $x$ is fixed, there is a price to pay in leveraging it to the ``best" feasibility guarantee in \eqref{main guarantee}. We discuss this next.

\subsection{Asymptotically Exact Coverage via $\chi^2$-Process Excursion}\label{sec:excursion}
The discussion so far presumes a fixed $x\in\Theta$. Recall in Section \ref{sec:good set} that, in data-driven DRO, a confidence region given by $\mathcal U$ guarantees $Z_0(x)\leq\max_{P\in\mathcal U}E_P[h(x;\xi)]$ with at least the same confidence level thanks to \eqref{good guarantee argument}. This guarantee holds regardless of a fixed $x$ or uniformly over all $x\in\Theta$. This is because the construction of such confidence regions is completely segregated from the expected value constraint of interest.
%
%Unlike the results in Section \ref{sec:EL} that only holds for a fixed $x\in\Theta$, the notion of asymptotically exact guarantee in \eqref{main guarantee}, however, notably requires the bound to hold for the whole feasible space $\Theta$. In standard data-driven DRO, a ``good" uncertainty set is constructed \emph{separately} from the function $h$, and because of this the statistical guarantee holds for all $x\in\Theta$ automatically, as revealed by the discussion in Section \ref{sec:good set}.
In contrast, the statistical performance of our empirical divergence ball is highly coupled with $h$, since $E_0[h(x;\xi)]$ can be viewed as the parameter we want to estimate in the EL method. Consequently, the reasoning for Theorem \ref{EL basic} only applies to situations where $x$ is fixed, and the empirical divergence ball constructed there is not big enough to guarantee \eqref{main guarantee}, which requires a bound simultaneous for all $x\in\Theta$.

The main result in this section is to explain and to show how, depending on the ``complexity" of $h$, one can suitably inflate the size of the ball to match a statistical performance close to \eqref{main guarantee}.

%To leverage to a simultaneous bound throughout all $\Theta$, we shall generalize the empirical likelihood machinery to the process level, and generalize the notion of empirical likelihood ratio
%
%Next we turn back to the stochastic program in \eqref{opt} where one needs to make decision $x$ in $Z(x)$.
We begin our discussion by imposing the following assumptions:
\begin{assumption}[Finite mean]
$Z_0(x)=E_0[h(x;\xi)]<\infty$ for all $x\in\Theta$.\label{finite mean}
\end{assumption}

\begin{assumption}[Non-degeneracy]
$\inf_{x\in\Theta}E_0|h(x;\xi)-Z_0(x)|>0$.\label{nondegeneracy}
\end{assumption}

\begin{assumption}[$L_2$-boundedness]
$E_0\sup_{x\in\Theta}|h(x;\xi)-Z_0(x)|^2<\infty$\label{L2}
\end{assumption}

\begin{assumption}[Function complexity]
The collection of functions
 \begin{equation}
 \mathcal H_\Theta=\{h(x;\cdot):\Xi\to\mathbb R|x\in\Theta\}\label{function class}
 \end{equation}
 is a $P_0$-Donsker class.\label{complexity}
\end{assumption}

The first three assumptions are mild moment conditions on the quantity $h(x;\xi)$. The last assumption, the so-called Donsker condition, means that the function class $\mathcal H_\Theta$ is ``simple" enough to allow the associated empirical process indexed by $\mathcal H_\Theta$ to converge weakly to a Brownian bridge (see Definition \ref{def:Donsker} in Appendix \ref{sec:EP}).

%\begin{assumption}
%The following conditions hold:
%\begin{enumerate}
% \item\emph{(Finite mean)} $Z_0(x)=E_0[h(x;\xi)]<\infty$ for all $x\in\Theta$.\label{finite mean}
% \item\emph{(Non-degeneracy)} $\inf_{x\in\Theta}E_0|h(x;\xi)-Z_0(x)|>0$.\label{nondegeneracy}
% \item\emph{($L_2$-boundedness)} $E_0\sup_{x\in\Theta}|h(x;\xi)-Z_0(x)|^2<\infty$\label{L2}
% \item\emph{(Function complexity)} The collection of functions
% \begin{equation}
% \mathcal H_\Theta=\{h(x;\cdot):\Xi\to\mathbb R|x\in\Theta\}\label{function class}
% \end{equation}
% is a $P_0$-Donsker class.\label{complexity}
% \end{enumerate}\label{assumptions}
%\end{assumption}

%For convenience, we denote $Z_n^*$ and $x_n^*$ as the optimal value and solution of the empirical DRO in \eqref{empirical opt}. The following provides a confidence bound guarantee for $x_n^*$:

The following theorem precisely describes the radius of the empirical divergence ball needed to attain the best guarantee in \eqref{main guarantee}:
\begin{theorem}[Optimal Calibration of Empirical Divergence Ball]
Let $\xi_1,\ldots,\xi_n\in\Xi$ be i.i.d. data under $P_0$. Suppose Assumptions \ref{finite mean}, \ref{nondegeneracy}, \ref{L2} and \ref{complexity} hold. Let $q_n$ be the $(1-\alpha)$-quantile of $\sup_{x\in\Theta}J_n(x)$, i.e.
\begin{equation}
P_{\bm\xi}\left(\sup_{x\in\Theta}J_n(x)\geq q_n\right)=\alpha\label{calibration}
\end{equation}
where $J_n(x)=G_n(x)^2$ and $G_n(\cdot)$ is a Gaussian process indexed by $\Theta$ that is centered, i.e. mean zero, with covariance
\begin{equation}
Cov(G_n(x_1),G_n(x_2))=\frac{\sum_{i=1}^n(h(x_1;\xi_i)-\hat h(x_1))(h(x_2;\xi_i)-\hat h(x_2))}{\sqrt{\sum_{i=1}^n(h(x_1;\xi_i)-\hat h(x_1))^2\sum_{i=1}^n(h(x_2;\xi_i)-\hat h(x_2))^2}}\label{cov}
\end{equation}
for any $x_1,x_2\in\Theta$, and $\hat h(x)=(1/n)\sum_{i=1}^nh(x;\xi_i)$ is the sample mean of $h(x;\xi_i)$'s. $P_{\bm\xi}$ denotes the probability conditional on the data $\xi_1,\ldots,\xi_n$.

We have
\begin{equation}
\lim_{n\to\infty}P(\underline Z_n^*(x)\leq Z_0(x)\leq\overline Z_n^*(x)\text{\ \ for all\ \ }x\in\Theta)=1-\alpha\label{asymptotic guarantee new}
\end{equation}
%where $Z_n^*$ and $x_n^*$ denotes the optimal value
where
\begin{align*}
\underline Z_n^*(x)&=\min_{\mathbf w\in\mathcal U_n(q_n/(2n))}\sum_{i=1}^nh(x;\xi_i)w_i\\
\overline Z_n^*(x)&=\max_{\mathbf w\in\mathcal U_n(q_n/(2n))}\sum_{i=1}^nh(x;\xi_i)w_i
\end{align*}
\label{main}
\end{theorem}

Note that, other than being a two-sided bound instead of one-sided, the guarantee \eqref{asymptotic guarantee new} is precisely \eqref{main guarantee}.

The process $J_n(\cdot)$, as the square of a Gaussian process, is known as a $\chi^2$-process (or $\chi^2$ random field; e.g., \cite{adler2009random}). Its covariance structure can be expressed explicitly in terms of the function $h$ and the data. The quantity $P_{\bm\xi}\left(\sup_{x\in\Theta}J_n(x)\geq u\right)$ is the excursion probability of $J_n(\cdot)$ above $u$. Note that we have ignored some subtle measurability issues in stating our result. To avoid unnecessary diversion, we will stay silent on measurability throughout the paper and refer the reader to \cite{van1996weak} for detailed treatments.

We observe some immediate connection of $\sup_xJ_n(x)$ to the $\chi^2_1$-distribution used in Theorem \ref{EL basic}. In addition to the fact that the marginal distribution of $J_n(x)$ at any $x$ is a $\chi^2_1$-distribution, we also have, by the Borell-TIS inequality (\cite{adler1990introduction}), that the asymptotic tail probability of $\sup_xJ_n(x)$ has the same exponential decay rate as that of $\chi^2_1$, i.e.
$$\frac{\log P\left(\sup_{x\in\Theta}G^2(x)\geq\nu\right)}{\log P(Y\geq\nu)}\to1$$
as $\nu\to\infty$, where $Y$ is a $\chi^2_1$ random variable. This suggests a relatively small overhead in using $q_n$ instead of $\chi^2_{1,1-\alpha}$ in calibrating the empirical ball when $\alpha$ is small.

Nevertheless, Theorem \ref{main} offers some insights beyond Theorem \ref{EL basic}. First, it requires the Donsker condition on the class $\mathcal H_\Theta$. %This condition requires the associated empirical process (REF) converges weakly to a Brownian bridge (the Gaussian process $G(\cdot)$). Intuitively, this holds when the function class $\mathcal H_\Theta$ is ``simple" in relation to the measure $P_0$. %We will clarify the precise definition of the Donsker class in Section \ref{sec:EP}. %, which is a measurement of the complexity of the function class $\mathcal H_\Theta$.
One sufficient condition of $P_0$-Donsker is:
\begin{lemma}
Suppose that $Z_0(x)=E_0[h(x;\xi)]<\infty$ and $Var_0(h(x;\xi))<\infty$ for all $x\in\Theta$. Also assume that there exists a random variable $M$ with $E_0M^2<\infty$ such that
$$|h(x_1;\xi)-h(x_2;\xi)|\leq M\|x_1-x_2\|_2$$
a.s. for all $x_1,x_2\in\Theta$. Then $\mathcal H_\Theta$ as defined in \eqref{function class} is $P_0$-Donsker.\label{Donsker sufficiency}
\end{lemma}
Lemma \ref{Donsker sufficiency} is a consequence of the Jain-Marcus Theorem (e.g., \cite{van1996weak}, Example 2.11.13). It is worth noting that the condition in Lemma \ref{Donsker sufficiency} is also a standard sufficient condition in guaranteeing the central limit convergence for SAA (\cite{shapiro2014lectures}, Theorem 5.7). This is not a coincidence, as the machinery behind Theorem \ref{main} involves an underpinning CLT, much like in the convergence analysis of SAA. %In this sense, suggests the empirical DRO is fundamentally closer to SAA than the standard data-driven DRO.

Secondly, even though $q_n\approx\chi^2_{1,1-\alpha}$ when $\alpha\approx0$, $q_n$ is strictly larger than $\chi^2_{1,1-\alpha}$ since $\sup_{x\in\Theta}J_n(x)$ stochastically dominates $\chi^2_1$ (unless degeneracy occurs). Thus the ball constructed in Theorem \ref{main} is always bigger than that in Theorem \ref{EL basic}. One way to estimate this inflation is by approximating the excursion probability of $\chi^2$-process using the theory of random geometry. We delegate this discussion to Section \ref{sec:Euler}. For now, we will delve into more details underlying Theorem \ref{main} and other properties of the empirical DRO.

\subsection{The Profile Nonparametric Likelihood Ratio Process and Other Properties of the Empirical DRO}\label{sec:PNLRP}
We explain briefly the machinery leading to Theorem \ref{main}, leaving the details to Appendix \ref{sec:proofs}. Our starting point is to define the profile nonparametric likelihood ratio in \eqref{EL process1} at the process level
\begin{equation}
\{R(x;Z):x\in\Theta\}\label{EL process}
\end{equation}
We call \eqref{EL process} the \emph{profile nonparametric likelihood ratio process} indexed by $x\in\Theta$. %, as
%\begin{equation}
%R(x;Z)=\max\left\{\prod_{i=1}^n(nw_i):\sum_{i=1}^nh(x;\xi_i)w_i=Z(x),\ \sum_{i=1}^nw_i=1,\ w_i\geq0\text{\ for\ }i=1,\ldots,n\right\}\label{EL process}
%\end{equation}
%and the profile likelihood ratio process indexed by $\Theta$ as the collection
%$$\{R(x):x\in\Theta\}$$
Denote the space
\begin{equation}
\ell^\infty(\Theta)=\left\{y:\Theta\to\mathbb R\Bigg|\|y\|_\Theta<\infty\right\}\label{topology}
\end{equation}
where we define $\|y\|_\Theta=\sup_{x\in\Theta}|y(x)|$ for any function $y:\Theta\to\mathbb R$. We have a convergence theorem for $R(x;Z)$ uniformly over $x\in\Theta$, in the following sense:

%\begin{lemma}
%Under Assumptions \ref{finite mean}, \ref{nondegeneracy}, \ref{L2} and \ref{complexity GC}, we have
%\begin{enumerate}
%\item \begin{equation}
%\left\|\frac{1}{n}\sum_{i=1}^n(h(\cdot;\xi_i)-Z_0(\cdot))^2-\sigma_0(\cdot)\right\|_\Theta\stackrel{a.s.}{\to}0\label{convergence var}
%\end{equation}
%where $\sigma_0(x)=Var_0(h(x;\xi))$.

\begin{theorem}[Limit Theorem of the Profile Nonparametric Likelihood Ratio Process]
%Assume that $Z_0(x):=E_0[h(x;\xi)]<\infty$ for all $x\in\Theta$, $\inf_{x\in\Theta}E_0|h(x;\xi)-Z_0(x)|\geq c>0$, $E\|h(\cdot;\xi)-Z_0(\cdot)\|_\Theta^2=E_0[\sup_{x\in\Theta}|h(x;\xi)-Z_0(x)|^2]<\infty$,
%and that $\mathcal H_\Theta=\{h(x;\cdot):x\in\Theta\}$ is a $P$-Donsker class.
Under Assumptions \ref{finite mean}, \ref{nondegeneracy}, \ref{L2} and \ref{complexity}, the profile likelihood ratio process defined in \eqref{EL process}
%$$R(x)=\sup\left\{\prod_{i=1}^n(nw_i):\sum_{i=1}^nh(x;\xi_i)w_i=Z(x),\ \sum_{i=1}^nw_i=1,\ w_i\geq0\text{\ for\ }i=1,\ldots,n\right\}$$
satisfies
$$-2\log R(\cdot;Z_0)\Rightarrow J(\cdot)\text{\ \ in\ \ }\ell^\infty(\Theta)$$
where $J(x)=G(x)^2$ and $G(\cdot)$ is a Gaussian process indexed by $x\in\Theta$ that has mean zero and covariance
$$Cov(G(x_1),G(x_2))=\frac{Cov_0(h(x_1;\xi),h(x_2;\xi))}{\sqrt{Var_0(h(x_1;\xi))Var_0(h(x_2;\xi))}}$$
for any $x_1,x_2\in\Theta$.%=\left\{x:\Theta\to\mathbb R\Bigg|\|x\|_{\Theta}=\sup_{\theta\in\Theta}|x(\theta)|<\infty\right\}$$
\label{EL process thm}
\end{theorem}

Theorem \ref{EL process thm} is the empirical-process generalization of Theorem \ref{EL basic}. It implies that $P\left(\sup_{x\in\Theta}\{-2\log R(x;Z_0)\}\leq q^*\right)\to1-\alpha$ for $q^*$ selected such that $P\left(\sup_{x\in\Theta}J(x)\leq q^*\right)=1-\alpha$. By a duality-type argument similar to that in Section \ref{sec:EL}, we have $-2\log R(x;Z_0)\leq q^*$ for all $x\in\Theta$, if and only if $\min_{\mathbf w\in\mathcal U_n(q^*/(2n))}\sum_{i=1}^nh(x;\xi_i)w_i\leq Z_0(x)\leq\max_{\mathbf w\in\mathcal U_n(q^*/(2n))}\sum_{i=1}^nh(x;\xi_i)w_i$ for all $x\in\Theta$, which implies $\lim_{n\to\infty}P\left(\min_{\mathbf w\in\mathcal U_n(q^*/(2n))}\sum_{i=1}^nh(x;\xi_i)w_i\leq Z_0(x)\leq\max_{\mathbf w\in\mathcal U_n(q^*/(2n))}\sum_{i=1}^nh(x;\xi_i)w_i\text{\ for all\ }x\in\Theta\right)=1-\alpha$ for the same choice of $q^*$. However, since $J$ relies on information about the unknown true distribution $P_0$, $q^*$ is unknown. The following result closes the gap by arguing that $J$ can be ``plugged-in" by $J_n$, and consequently $q^*$ by $q_n$ as depicted in Theorem \ref{main}:
\begin{lemma}
Under Assumptions \ref{finite mean}, \ref{nondegeneracy}, \ref{L2} and \ref{complexity}, conditional on almost every data realization $(P_n:n\geq1)$,
$$G_n(\cdot)\Rightarrow G(\cdot)\text{\ \ in\ \ }\ell^\infty(\Theta)$$
where $G_n(\cdot)$ and $G(\cdot)$ are defined in Theorems \ref{main} and \ref{EL process thm} respectively.\label{sample Gaussian}
\end{lemma}

Theorem \ref{main} can then be proved by combining Theorem \ref{EL process thm} and Lemma \ref{sample Gaussian}. %, and subsequently Theorem \ref{main} in Section \ref{sec:proofs}.
Moreover, consistency of the empirical DRO also holds uniformly over $x\in\Theta$:
\begin{theorem}[Uniform Strong Consistency]
Under Assumptions \ref{finite mean}, \ref{nondegeneracy}, \ref{L2} and \ref{complexity},
\begin{align*}
\sup_{x\in\Theta}|\underline Z_n^*(x)-Z_0(x)|\stackrel{a.s.}{\to}0\\
\sup_{x\in\Theta}|\overline Z_n^*(x)-Z_0(x)|\stackrel{a.s.}{\to}0
\end{align*}
as $n\to\infty$.
\label{consistency}
\end{theorem}

Lastly, the following theorem highlights that the width of the confidence band $[\underline Z_n^*(x),\overline Z_n^*(x)]$ varies with the standard deviation at each $x$:
\begin{theorem}[Pertaining to the Variability at Each Decision Point]
Suppose Assumptions \ref{finite mean}, \ref{nondegeneracy}, \ref{L2} and \ref{complexity} hold. Additionally, suppose that $h(\cdot;\cdot)$ is bounded. Then
$$\underline Z_n^*(x)=\hat h(x)-\sqrt{q_n}\frac{\hat\sigma(x)}{\sqrt n}+O\left(\frac{1}{n}\right)$$
$$\overline Z_n^*(x)=\hat h(x)+\sqrt{q_n}\frac{\hat\sigma(x)}{\sqrt n}+O\left(\frac{1}{n}\right)$$
uniformly over $x\in\Theta$ a.s.. Here $\hat h(x)=\frac{1}{n}\sum_{i=1}^nh(x;\xi_i)$ is the sample mean, $\hat\sigma^2(x)=\frac{1}{n}\sum_{i=1}^n(h(x;\xi_i)-\hat h(x))^2$ is the sample variance at each $x$, and $q_n$ is defined in Theorem \ref{main}.\label{asymptotic equivalence}
\end{theorem}
Theorem \ref{asymptotic equivalence} gives rise to \eqref{mimic CI}. In particular, $\sqrt{q_n}$ is analogous to the critical value in a confidence band. In summary, Theorems \ref{main} and \ref{asymptotic equivalence} show that our empirical divergence ball $\mathcal U_n(q_n/(2n))$, calibrated via the quantile of $\chi^2$-process excursion $q_n$, satisfies our benchmark guarantees \eqref{main guarantee} and \eqref{mimic CI}, except that it provides a two-sided bound instead of one-sided. The difference of two- versus one-sided bound is the reason we have claimed ``close to" the best in Section \ref{sec:good set}.

%This theorem shows that $[\underline Z_n^*(x),\overline Z_n^*(x)]$ indeed mimics a simultaneous confidence band built from the standard (functional) CLT, hence concluding our claim in Section \ref{sec:good set}. This result, by now, should not be surprising since we have shown asymptotically exact coverage guarantee when fixing $x$, which implies, at least intuitively, that empirical DRO gives rise to a structure that is asymptotically the same as CLT. Note that Theorem \ref{asymptotic equivalence} requires the additional assumption of a bounded $h$ relative to all our previous results.

\subsection{Approximating the Quantile of $\chi^2$-process Excursion}\label{sec:Euler}
We discuss how to estimate $q_n$ in Theorem \ref{main}. One approach is to approximate the excursion probability of $\chi^2$-process by the mean Euler characteristic approximation (e.g., \cite{adler2009random}, Theorem 13.4.1 and Section 15.10.2, and \cite{adler2011topological}, Theorem 4.8.1):
%using the mean Euler characteristic approximation (REF)
\begin{equation}
P\left(\sup_{x\in\Theta}J_n(x)\geq u\right)\approx\sum_{j=0}^m(2\pi)^{-j/2}\mathcal L_j(\Theta)\mathcal M_j(u)\label{Euler}
\end{equation}
%E[\varphi(A_D(G_n,\Theta))]=\sum_{j=0}^m(2\pi)^{-j/2}\mathcal L_j(\Theta)\mathcal M_j(D)\label{Euler}
%\end{equation}
Here $m$ is the dimension of the decision space $\Theta\subset\mathbb R^m$. The coefficients $\mathcal L_j(\Theta)$ on the RHS of \eqref{Euler} are known as the Lipschitz-Killing curvatures of the domain $\Theta$, which measure the ``intrinsic volumes" of the domain $\Theta$ using the Riemannian metric induced by the Gaussian process $G_n$ (\cite{adler2009random}, equation (12.2.2)). In particular, the highest-dimensional coefficient is given by
$$\mathcal L_m(\Theta)=\int_\Theta det(\Lambda(x))^{1/2}dx$$
(\cite{adler2009random}, equation (12.2.22), and \cite{adler2011topological}, equation (5.4.1)) where $\Lambda(x)=(\Lambda_{ij}(x))_{i,j=1,\ldots,m}\in\mathbb R^{m\times m}$, and %is the second order spectral moment of $G_n$, namely
\begin{equation}
\Lambda_{ij}(x)=Cov\left(\frac{\partial G_n(x)}{\partial x_i},\frac{\partial G_n(x)}{\partial x_j}\right)=\frac{\partial^2}{\partial y_i\partial z_j}Cov(G_n(y),G_n(z))\Big|_{y=x,z=x}\label{general moment}
\end{equation}
for differentiable $G_n$ (in the $L^2$ sense), with $x_i$ and $x_j$ the $i$ and $j$-th components of $x$. %This is equal to
%\begin{equation}
%\label
%\end{equation}
%Furthermore, when $\Theta$ is one-dimensional, and $G_n$ is centered and has unit variance, we have
%\begin{equation}
%\Lambda(x)=-\frac{\partial^2}{\partial y^2}Cov(G_n(y),G_n(x))\Big|_{y=x}\label{one-dim moment}
%\end{equation}
Thus \eqref{general moment} can be evaluated by differentiating \eqref{cov}. Lower-dimensional coefficients can be evaluated by integration over lower-dimensional surfaces of $\Theta$, and $\mathcal L_0(\Theta)=1$.

On the other hand, the quantities $\mathcal M_j(u)$'s are the Gaussian Minkowski functionals for the excursion set, independent of $\Theta$ and $h$, and are given by
$$\mathcal M_j(u)=(-1)^j\frac{d^j}{dy^j}P(Y\geq y)\Big|_{y=\sqrt u}$$
where $Y$ is the square root of a $\chi^2_1$ random variable. Thus, for instance, $\mathcal M_0(u)=P(\chi^2_1\geq u)$, and $\mathcal M_1(u)=2\phi(\sqrt u)$ where $\phi(\cdot)$ is the standard normal density.

%The approximation \eqref{Euler} is known as the mean Euler characteristic approximation, developed in the theory of random geometry (e.g., \cite{adler2009random}, Theorem 13.4.1 and Section 15.10.2, and \cite{adler2011topological} Theorem 4.8.1). %see, e.g., \cite{adler2009random}, Chapters 12--15).
\eqref{Euler} is a very accurate approximation for $P(\sup_{x\in\Theta}J_n(x)\geq u)$ in the sense
$$\left|P\left(\sup_{x\in\Theta}J_n(x)\geq u\right)-\sum_{j=0}^m(2\pi)^{-j/2}\mathcal L_j(\Theta)\mathcal M_j(u)\right|\leq Ce^{-\beta u/2}$$
where $C>0$ and $\beta>1$ (\cite{atw16}, Section 5.3.2). In other words, the approximation error is exponentially smaller than all the terms in \eqref{Euler} as $u$ increases. In practice, however, the formula for $\mathcal L_j(\Theta)$ could get increasingly complex as $j$ decreases, in which case only the first and the second highest order coefficients of $\mathcal L_j(\Theta)$ are used.

In light of the above, an accurate approximation of $q_n$ can be found by solving the root of
$$\sum_{j=0}^m(2\pi)^{-j/2}\mathcal L_j(\Theta)\mathcal M_j(u)=\alpha$$

%Here $D=\{y\in\mathbb R:y^2\geq u\}$ for some $u>0$, and $A_D(G_n,\Theta)=\{x\in\Theta:G_n(x)\in D\}$ is the excursion set for $G_n$ associated with $D$. Equivalently, $A_D(G_n,\Theta)=\{x\in\Theta:J_n(x)\geq u\}$ is the excursion set of $J_n$ above the level $u$. $\varphi:\Theta\to\mathbb R$ defines the so-called \emph{Euler characteristic} of an excursion set, a well-studied fundamental geometric property (REF).

As an explicit illustration, when $m=1$, and $h$ is twice differentiable almost everywhere, we have
$$P\left(\sup_{x\in\Theta}J_n(x)\geq u\right)=P(\chi^2_1\geq u)+\int_\Theta \sqrt{\frac{\partial^2}{\partial y\partial z}Cov(G_n(y),G_n(z))\Big|_{y=x,z=x}}dx\ \frac{e^{-u/2}}{\pi}+O(e^{-\beta u/2})$$
for some $\beta>1$. An approximate $q_n$ can then be found by solving the root of
\begin{equation}
P(\chi^2_1\geq u)+\int_\Theta \sqrt{\frac{\partial^2}{\partial y\partial z}Cov(G_n(y),G_n(z))\Big|_{y=x,z=x}}dx\ \frac{e^{-u/2}}{\pi}=\alpha\label{calibrate inflated}
\end{equation}

\section{Numerical Illustrations for the Empirical DRO}\label{sec:numerics}
This section shows some numerical results on the statistical performance of empirical DRO. We use the newsvendor loss function in \eqref{h example} as our $h$. We repeat $1,000$ times:
\begin{enumerate}
\item Simulate $n$ i.i.d. data $\xi_1,\ldots,\xi_n$ from the $k$-discretized $Exp(1/20)$.
\item Estimate $q_n$ using \eqref{calibrate inflated}, and compute $\underline Z_n^*(x)=\min_{\mathbf w\in\mathcal U_n(q_n/(2n))}\sum_{i=1}^nh(x;\xi_i)w_i$ and $\overline Z_n^*(x)=\max_{\mathbf w\in\mathcal U_n(q_n/(2n))}\sum_{i=1}^nh(x;\xi_i)w_i$, with $\alpha$ set to be $0.05$.
\item Output
$$I\left(\underline Z_n^*(x)\leq Z_0(x)\leq\overline Z_n^*(x)\text{\ \ for\ }x=\frac{50j}{20},\ j=1,\ldots,20\right)$$
and
$$I\left(Z_0(x)\leq\overline Z_n^*(x)\text{\ \ for\ }x=\frac{50j}{20},\ j=1,\ldots,20\right)$$
where $Z_0(x)$ is the true function of interest that is calculable in closed-form.
\end{enumerate}
We then output the point estimates and the $95\%$ CIs of the two- and one-sided coverage probabilities, i.e. $P\left(\underline Z_n^*(x)\leq Z_0(x)\leq\overline Z_n^*(x)\text{\ \ for\ }x=\frac{50j}{20},\ j=1,\ldots,20\right)$ and $P\left(Z_0(x)\leq\overline Z_n^*(x)\text{\ \ for\ }x=\frac{50j}{20},\ j=1,\ldots,20\right)$, from the $1,000$ replications. These probabilities serve as proxies for the probabilities $P\left(\underline Z_n^*(x)\leq Z_0(x)\leq\overline Z_n^*(x)\text{\ \ for\ }x\in\Theta\right)$ and $P\left(Z_0(x)\leq\overline Z_n^*(x)\text{\ \ for\ }x\in\Theta\right)$ respectively.

We first set $\xi$ as a $k$-discretized $Exp(1/20)$ as in \eqref{exp specification}. For comparison, we also repeat the above experiment but using $\chi^2_{k-1,0.95}$ and $\chi^2_{1,0.95}$ in place of $q_n$. Table \ref{table:comparison discrete n} shows the results of two-sided coverage probabilities as we vary the sample size from $n=20$ to $80$. The coverage probabilities appear to be stable already starting at $n=20$. As we can see, the coverages using the $\chi^2_{k-1,0.95}$ calibration (Table \ref{table:DRO simultaneous n}) are around $99\%$, much higher than $95\%$, as $k-1$ is over-determining the number of parameters we want to estimate from the EL perspective. The coverage probabilities using the $\chi^2_{1,0.95}$ calibration (Table \ref{table:CLT simultaneous n}), on the other hand, are in the range $86\%$ to $87\%$, significantly lower than $95\%$, since it does not account for simultaneous estimation errors. Lastly, the coverage probabilities using the $\chi^2$-process excursion (Table \ref{table:empirical simultaneous n}) are very close to $95\%$ in all cases, thus confirming the superiority of our approach.

\begin{table}[!htb]
  \begin{subtable}{0.33\textwidth}
    \centering
    {\footnotesize\begin{tabular}{c|c|c}
    $n$&2-sided&$95\%$ C.I. of\\
    &Cover. Prob.&Cover. Prob.\\
    \hline
   20&98.3\%&(98.0\%, 98.6\%)\\
30&98.8\%&(98.7\%, 100.0\%)\\
40&98.9\%&(98.8\%, 100.0\%)\\
50&98.8\%&(98.7\%, 100.0\%)\\
60&98.9\%&(98.8\%, 100.0\%)\\
80&98.8\%&(98.6\%, 98.9\%)
    \end{tabular}}
    \subcaption{$\chi^2_{k-1,0.95}/(2n)$}
     %     \captionof{table}{Coverage probabilities for distributions of different support sizes using data-driven DRO}
      \label{table:DRO simultaneous n}
  \end{subtable}
  %\hspace{.3cm}
  \begin{subtable}{0.33\textwidth}
    \centering
    {\footnotesize\begin{tabular}{c|c|c}
    $n$&2-sided&$95\%$ C.I. of\\
    &Cover. Prob.&Cover. Prob.\\
    \hline
20&85.6\%&(85.0\%, 86.2\%)\\
30&85.9\%&(85.4\%, 86.5\%)\\
40&86.6\%&(86.1\%, 87.1\%)\\
50&86.6\%&(86.1\%, 87.1\%)\\
60&86.1\%&(85.5\%, 86.7\%)\\
80&86.8\%&(86.3\%, 87.3\%)
    \end{tabular}}
    \subcaption{$\chi^2_{1,0.95}/(2n)$}
     %     \captionof{table}{Coverage probabilities for distributions of different support sizes using standard CLT}
      \label{table:CLT simultaneous n}
    \end{subtable}
  \begin{subtable}{0.33\textwidth}
    \centering
    {\footnotesize\begin{tabular}{c|c|c}
    $n$&2-sided&$95\%$ C.I. of\\
    &Cover. Prob.&Cover. Prob.\\
    \hline
20&94.4\%&(93.8\%, 95.0\%)\\
30&94.6\%&(94.0\%, 95.2\%)\\
40&94.7\%&(94.2\%, 95.3\%)\\
50&94.7\%&(94.2\%, 95.3\%)\\
60&94.4\%&(93.8\%, 94.9\%)\\
80&95.0\%&(94.5\%, 95.5\%)
    \end{tabular}}
    \subcaption{Approximate $95\%$-quantile of $\sup_xJ_n(x)$}
     %     \captionof{table}{Coverage probabilities for distributions of different support sizes using empirical DRO}
      \label{table:empirical simultaneous n}
       \end{subtable}
       \caption{Two-sided coverage probabilities for different Burg-divergence ball sizes and sample sizes for a discrete distribution with $k=5$}
\label{table:comparison discrete n}
  \end{table}

Next, Table \ref{table:comparison discrete n 1-sided} shows the results for one-sided coverage instead of two-sided. These one-sided coverage probabilities are slightly higher than the two-sided counterparts as the coverage condition is now more relaxed. Nonetheless, the magnitudes of these changes are very small compared to the effects brought by the choice of calibration methods. In particular, using $\chi^2_{k-1,0.95}$ appears to be severely over-covering at about $99\%$ to $100\%$, while using $\chi^2_{1,0.95}$ gives under-coverage at about $89\%$ to $91\%$. Using the $\chi^2$-process excursion shows $95\%$ to $96\%$ coverage performances, thus significantly better than the other two methods. These show that, even though our statistical guarantees in Theorem \ref{main} are two-sided, the loss of inaccuracy for one-sided coverage is very minor compared to the improvement in the calibration method used. This experimentally justifies our claim of ``close to" the best at the end of Section \ref{sec:PNLRP}.

\begin{table}[!htb]
  \begin{subtable}{0.33\textwidth}
    \centering
    {\footnotesize\begin{tabular}{c|c|c}
    $n$&1-sided&$95\%$ C.I. of\\
    &Cover. Prob.&Cover. Prob.\\
    \hline
     20&98.6\%&(98.1\%,	99.2\%)\\
30&99.6\%&(99.4\%, 100.0\%)\\
40&99.7\%&(99.6\%, 100.0\%)\\
50&99.6\%&(99.4\%, 100.0\%)\\
60&99.7\%&(99.6\%, 100.0\%)\\
80&99.5\%&(99.3\%, 99.8\%)
    \end{tabular}}
    \subcaption{$\chi^2_{k-1,0.95}/(2n)$}
     %     \captionof{table}{Coverage probabilities for distributions of different support sizes using data-driven DRO}
      \label{table:DRO simultaneous n 1-sided}
  \end{subtable}
  %\hspace{.3cm}
  \begin{subtable}{0.33\textwidth}
    \centering
    {\footnotesize\begin{tabular}{c|c|c}
    $n$&1-sided&$95\%$ C.I. of\\
    &Cover. Prob.&Cover. Prob.\\
    \hline
           20&88.5\%&(87.2\%, 89.8\%)\\
30&89.2\%&(88.0\%, 90.4\%)\\
40&90.5\%&(89.4\%, 91.6\%)\\
50&90.5\%&(89.4\%, 91.6\%)\\
60&89.5\%&(88.3\%, 90.7\%)\\
80&90.9\%&(89.9\%, 91.9\%)
    \end{tabular}}
    \subcaption{$\chi^2_{1,0.95}/(2n)$}
     %     \captionof{table}{Coverage probabilities for distributions of different support sizes using standard CLT}
      \label{table:CLT simultaneous n 1-sided}
    \end{subtable}
  \begin{subtable}{0.33\textwidth}
    \centering
    {\footnotesize\begin{tabular}{c|c|c}
    $n$&1-sided&$95\%$ C.I. of\\
    &Cover. Prob.&Cover. Prob.\\
    \hline
       20&94.8\%&(93.7\%, 96.0\%)\\
30&95.1\%&(94.0\%, 96.3\%)\\
40&95.4\%&(94.3\%, 96.5\%)\\
50&95.4\%&(94.3\%, 96.5\%)\\
60&94.7\%&(93.5\%, 95.8\%)\\
80&96.0\%&(95.0\%, 97.0\%)
    \end{tabular}}
    \subcaption{Approximate $95\%$-quantile of $\sup_xJ_n(x)$}
     %     \captionof{table}{Coverage probabilities for distributions of different support sizes using empirical DRO}
      \label{table:empirical simultaneous n 1-sided}
       \end{subtable}
       \caption{One-sided coverage probabilities for different Burg-divergence ball sizes and sample sizes for a discrete distribution with $k=5$}
\label{table:comparison discrete n 1-sided}
  \end{table}

Finally, we repeat the experiments using the continuous distribution $Exp(1/20)$. We compare the use of $\chi^2_{1,0.95}$ with the $\chi^2$-process excursion (there is no notion of $k$ in this case). Table \ref{table:comparison cont} shows that the two-sided coverages using $\chi^2_{1,0.95}$ are under-covering at between $82\%$ and $85\%$. The $\chi^2$-process excursion gives about $93\%$ at $n=20$ and converges to close to $95\%$ at $n=80$. Thus, similar to the discrete case, the calibration using $\chi^2$-process excursion gives significantly more accurate two-sided coverages than using $\chi^2_{1,0.95}$. Table \ref{table:comparison cont 1-sided} draws similar conclusion for one-sided coverages. For $\chi^2_{1,0.95}$, the coverage probability is about $84\%$ at $n=20$ and $90\%$ at $n=80$, therefore severely under-covering. On the other hand, $\chi^2$-process excursion gives $94\%$ to $96\%$ coverages among all the $n$'s. This once again shows the insignificance of one- versus two-sided coverage compared to the improvement in the choice of calibration method. In overall, our proposed scheme of using $\chi^2$-process excursion gives much more accurate coverages than using $\chi^2_{1,0.95}$.

\begin{table}[!htb]
  \begin{subtable}{0.5\textwidth}
    \centering
    {\footnotesize\begin{tabular}{c|c|c}
    $n$&2-sided&$95\%$ C.I. of\\
    &Cover. Prob.&Cover. Prob.\\
    \hline
  20&82.4\%&(81.3\%, 83.5\%)\\
30&82.8\%&(81.7\%, 83.8\%)\\
40&83.5\%&(82.5\%, 84.5\%)\\
50&84.3\%&(83.3\%, 85.2\%)\\
60&84.2\%&(83.2\%, 85.2\%)\\
80&85.1\%&(84.2\%, 86.1\%)
   \end{tabular}}
   \subcaption{$\chi^2_{1,0.95}/(2n)$}
      %     \captionof{table}{Coverage probabilities for distributions of different support sizes using data-driven DRO}
      \label{table:DRO simultaneous cont}
  \end{subtable}
  %\hspace{.3cm}
  \begin{subtable}{0.5\textwidth}
    \centering
    {\footnotesize\begin{tabular}{c|c|c}
    $n$&2-sided&$95\%$ C.I. of\\
    &Cover. Prob.&Cover. Prob.\\
    \hline
    20&92.8\%&(91.1\%, 94.6\%)\\
30&93.6\%&(92.3\%, 94.9\%)\\
40&94.4\%&(93.1\%, 95.8\%)\\
50&94.4\%&(93.1\%, 95.8\%)\\
60&95.7\%&(94.6\%, 96.9\%)\\
80&95.3\%&(94.1\%, 96.5\%)
 \end{tabular}}
 \subcaption{Approximate $95\%$-quantile of $\sup_xJ_n(x)$}
         %     \captionof{table}{Coverage probabilities for distributions of different support sizes using standard CLT}
      \label{table:empirical simultaneous cont}
    \end{subtable}
    \caption{Two-sided coverage probabilities for different Burg-divergence ball sizes and sample sizes for a continuous distribution}
\label{table:comparison cont}
   \end{table}

\begin{table}[!htb]
  \begin{subtable}{0.5\textwidth}
    \centering
    {\footnotesize\begin{tabular}{c|c|c}
    $n$&1-sided&$95\%$ C.I. of\\
    &Cover. Prob.&Cover. Prob.\\
    \hline
      20&83.9\%&(81.7\%, 86.1\%)\\
30&84.6\%&(82.4\%, 86.8\%)\\
40&86.1\%&(84.0\%, 88.2\%)\\
50&87.7\%&(85.7\%, 89.7\%)\\
60&87.5\%&(85.4\%, 89.6\%)\\
80&89.5\%&(87.6\%, 91.4\%)
   \end{tabular}}
   \subcaption{$\chi^2_{1,0.95}/(2n)$}
      %     \captionof{table}{Coverage probabilities for distributions of different support sizes using data-driven DRO}
      \label{table:DRO simultaneous cont 1-sided}
  \end{subtable}
  %\hspace{.3cm}
  \begin{subtable}{0.5\textwidth}
    \centering
    {\footnotesize\begin{tabular}{c|c|c}
    $n$&1-sided&$95\%$ C.I. of\\
    &Cover. Prob.&Cover. Prob.\\
    \hline
    20&93.6\%&(91.9\%, 95.3\%)\\
30&94.4\%&(93.1\%, 95.7\%)\\
40&95.2\%&(93.9\%, 96.5\%)\\
50&95.2\%&(93.9\%, 96.5\%)\\
60&96.5\%&(95.4\%, 97.6\%)\\
80&96.1\%&(94.9\%, 97.3\%)
 \end{tabular}}
 \subcaption{Approximate $95\%$-quantile of $\sup_xJ_n(x)$}
         %     \captionof{table}{Coverage probabilities for distributions of different support sizes using standard CLT}
      \label{table:empirical simultaneous cont 1-sided}
    \end{subtable}
    \caption{One-sided coverage probabilities for different Burg-divergence ball sizes and sample sizes for a continuous distribution}
\label{table:comparison cont 1-sided}
   \end{table}

%Lastly, Figures ? and ? show the values of $\underline Z_n^*(x)$ and $\overline Z_n^*(x)$ as a function of $x\in[0,50]$. We can see that the width of the interval varies according to the value of $x$.

\section{Conclusion}\label{sec:conlusion}
We have motivated and investigated the construction of tractable uncertainty sets that can recover the feasibility guarantees on par with the implications of CLT. We have shown that the empirical Burg-entropy divergence balls are capable of achieving such guarantees. We have also shown, intriguingly, that these balls are invalid confidence regions in the standard framework of data-driven DRO, and can have low or zero coverages on the true underlying distributions. Rather, we have explained their statistical performances via linking the resulting DRO with empirical likelihood. This link allows us to derive the optimal sizes of these balls, using the quantiles of $\chi^2$-process excursion. Such a calibration approach also bypasses some documented difficulties in using divergence balls in the data-driven DRO literature. Future work includes further developments of the theory and calibration methods to incorporate optimization objectives and more general constraints.

\ACKNOWLEDGMENT{The author gratefully acknowledges support from the National Science Foundation under grants CMMI-1400391/1542020 and CMMI-1436247/1523453. He also thanks Zhiyuan Huang and Yanzhe Jin for assisting with the numerical experiments.}

\bibliographystyle{informs2014} % outcomment this and next line in Case 1
\bibliography{bibliography} % if more than one, comma separated

\begin{thebibliography}{34}
\providecommand{\natexlab}[1]{#1}
\providecommand{\url}[1]{\texttt{#1}}
\providecommand{\urlprefix}{URL }

\bibitem[{Adler \protect\BIBand{} Taylor(2011)}]{adler2011topological}
Adler R, Taylor JE (2011) \emph{Topological Complexity of Smooth Random
  Functions: {\'E}cole D'{\'E}t{\'e} de Probabilit{\'e}s de Saint-Flour
  XXXIX-2009} (Springer Science \& Business Media).

\bibitem[{Adler(1990)}]{adler1990introduction}
Adler RJ (1990) \emph{An Introduction to Continuity, Extrema, and Related
  Topics for General {G}aussian Processes} (Institute of Mathematical
  Statistics, Lecture Notes-Monograph Series).

\bibitem[{Adler \protect\BIBand{} Taylor(2009)}]{adler2009random}
Adler RJ, Taylor JE (2009) \emph{Random Fields and Geometry} (Springer Science
  \& Business Media).

\bibitem[{Adler et~al.(in preparation)Adler, Taylor, \protect\BIBand{}
  Worsley}]{atw16}
Adler RJ, Taylor JE, Worsley KJ (in preparation) \emph{Applications of Random
  Fields and Geometry: Foundations and Case Studies} (available at
  http://webee.technion.ac.il/people/adler/hrf.pdf).

\bibitem[{Agresti \protect\BIBand{} Kateri(2011)}]{agresti2011categorical}
Agresti A, Kateri M (2011) \emph{Categorical Data Analysis} (Springer).

\bibitem[{Atlason et~al.(2004)Atlason, Epelman, \protect\BIBand{}
  Henderson}]{atlason2004call}
Atlason J, Epelman MA, Henderson SG (2004) Call center staffing with simulation
  and cutting plane methods. \emph{Annals of Operations Research}
  127(1-4):333--358.

\bibitem[{Ben-Tal et~al.(2013)Ben-Tal, Den~Hertog, De~Waegenaere, Melenberg,
  \protect\BIBand{} Rennen}]{ben2013robust}
Ben-Tal A, Den~Hertog D, De~Waegenaere A, Melenberg B, Rennen G (2013) Robust
  solutions of optimization problems affected by uncertain probabilities.
  \emph{Management Science} 59(2):341--357.

\bibitem[{Bertsimas et~al.(2014)Bertsimas, Gupta, \protect\BIBand{}
  Kallus}]{bertsimas2014robust}
Bertsimas D, Gupta V, Kallus N (2014) Robust {SAA}. \emph{arXiv preprint
  arXiv:1408.4445} .

\bibitem[{Cox \protect\BIBand{} Hinkley(1979)}]{cox1979theoretical}
Cox DR, Hinkley DV (1979) \emph{Theoretical Statistics} (CRC Press).

\bibitem[{Delage \protect\BIBand{} Ye(2010)}]{delage2010distributionally}
Delage E, Ye Y (2010) Distributionally robust optimization under moment
  uncertainty with application to data-driven problems. \emph{Operations
  research} 58(3):595--612.

\bibitem[{Esfahani \protect\BIBand{} Kuhn(2015)}]{esfahani2015data}
Esfahani PM, Kuhn D (2015) Data-driven distributionally robust optimization
  using the {W}asserstein metric: Performance guarantees and tractable
  reformulations. \emph{arXiv preprint arXiv:1505.05116} .

\bibitem[{F{\'a}bi{\'a}n(2008)}]{fabian2008handling}
F{\'a}bi{\'a}n CI (2008) Handling {CV}a{R} objectives and constraints in
  two-stage stochastic models. \emph{European Journal of Operational Research}
  191(3):888--911.

\bibitem[{Gabrel et~al.(2014)Gabrel, Murat, \protect\BIBand{}
  Thiele}]{gabrel2014recent}
Gabrel V, Murat C, Thiele A (2014) Recent advances in robust optimization: An
  overview. \emph{European Journal of Operational Research} 235(3):471--483.

\bibitem[{Gao \protect\BIBand{} Kleywegt(2016)}]{gao2016distributionally}
Gao R, Kleywegt AJ (2016) Distributionally robust stochastic optimization with
  {W}asserstein distance. \emph{arXiv preprint arXiv:1604.02199} .

\bibitem[{Goh \protect\BIBand{} Sim(2010)}]{goh2010distributionally}
Goh J, Sim M (2010) Distributionally robust optimization and its tractable
  approximations. \emph{Operations research} 58(4-part-1):902--917.

\bibitem[{Gupta(2015)}]{gupta2015near}
Gupta V (2015) Near-optimal ambiguity sets for distributionally robust
  optimization. \emph{Preprint} .

\bibitem[{Jiang \protect\BIBand{} Guan(2012)}]{jiang2012data}
Jiang R, Guan Y (2012) Data-driven chance constrained stochastic program.
  \emph{Mathematical Programming} 1--37.

\bibitem[{Kleywegt et~al.(2002)Kleywegt, Shapiro, \protect\BIBand{} Homem-de
  Mello}]{kleywegt2002sample}
Kleywegt AJ, Shapiro A, Homem-de Mello T (2002) The sample average
  approximation method for stochastic discrete optimization. \emph{SIAM Journal
  on Optimization} 12(2):479--502.

\bibitem[{Krokhmal et~al.(2002)Krokhmal, J~Palmquist, \protect\BIBand{}
  Uryasev}]{palmquist1999portfolio}
Krokhmal P, J~Palmquist J, Uryasev S (2002) Portfolio optimization with
  conditional value-at-risk objective and constraints. \emph{Journal of Risk}
  4:43--68.

\bibitem[{Kullback \protect\BIBand{} Leibler(1951)}]{kullback1951information}
Kullback S, Leibler RA (1951) On information and sufficiency. \emph{The Annals
  of Mathematical Statistics} 22(1):79--86.

\bibitem[{Lam \protect\BIBand{} Zhou(2015)}]{lam2015quantifying}
Lam H, Zhou E (2015) Quantifying uncertainty in sample average approximation.
  \emph{Proceedings of the 2015 Winter Simulation Conference}, 3846--3857 (IEEE
  Press).

\bibitem[{Moon \protect\BIBand{} Hero(2014)}]{moon2014multivariate}
Moon K, Hero A (2014) Multivariate f-divergence estimation with confidence.
  \emph{Advances in Neural Information Processing Systems}, 2420--2428.

\bibitem[{Nagaraj \protect\BIBand{}
  Pasupathy(2014)}]{nagaraj2014stochastically}
Nagaraj K, Pasupathy R (2014) Stochastically constrained simulation
  optimization on integer-ordered spaces: The cg{R}-{SPLINE} algorithm.
  \emph{Technical Report. Purdue University, West Lafayette, IN} .

\bibitem[{Nguyen et~al.(2007)Nguyen, Wainwright, \protect\BIBand{}
  Jordan}]{nguyen2007estimating}
Nguyen X, Wainwright MJ, Jordan MI (2007) Estimating divergence functionals and
  the likelihood ratio by penalized convex risk minimization. \emph{NIPS},
  1089--1096.

\bibitem[{Owen(1988)}]{owen1988empirical}
Owen AB (1988) Empirical likelihood ratio confidence intervals for a single
  functional. \emph{Biometrika} 75(2):237--249.

\bibitem[{Owen(2001)}]{owen2001empirical}
Owen AB (2001) \emph{Empirical likelihood} (CRC press).

\bibitem[{P{\'a}l et~al.(2010)P{\'a}l, P{\'o}czos, \protect\BIBand{}
  Szepesv{\'a}ri}]{pal2010estimation}
P{\'a}l D, P{\'o}czos B, Szepesv{\'a}ri C (2010) Estimation of {R}{\'e}nyi
  entropy and mutual information based on generalized nearest-neighbor graphs.
  \emph{Advances in Neural Information Processing Systems}, 1849--1857.

\bibitem[{Pardo(2005)}]{pardo2005statistical}
Pardo L (2005) \emph{Statistical Inference Based on Divergence Measures} (CRC
  Press).

\bibitem[{Shapiro et~al.(2014)Shapiro, Dentcheva, \protect\BIBand{}
  Ruszczy{\'n}ski}]{shapiro2014lectures}
Shapiro A, Dentcheva D, Ruszczy{\'n}ski A (2014) \emph{Lectures on Stochastic
  Programming: Modeling and Theory}, volume~16 (SIAM).

\bibitem[{Van Der~Vaart \protect\BIBand{} Wellner(2000)}]{van2000preservation}
Van Der~Vaart A, Wellner JA (2000) Preservation theorems for
  {G}livenko-{C}antelli and uniform {G}livenko-{C}antelli classes. \emph{High
  dimensional probability II}, 115--133 (Springer).

\bibitem[{Van Der~Vaart \protect\BIBand{} Wellner(1996)}]{van1996weak}
Van Der~Vaart AW, Wellner JA (1996) \emph{Weak Convergence and Empirical
  Processes: With Applications in Statistics} (Springer).

\bibitem[{Wang \protect\BIBand{} Ahmed(2008)}]{wang2008sample}
Wang W, Ahmed S (2008) Sample average approximation of expected value
  constrained stochastic programs. \emph{Operations Research Letters}
  36(5):515--519.

\bibitem[{Wang et~al.(2015)Wang, Glynn, \protect\BIBand{} Ye}]{Wang2015}
Wang Z, Glynn PW, Ye Y (2015) Likelihood robust optimization for data-driven
  problems. \emph{Computational Management Science} 13(2):241--261.

\bibitem[{Wiesemann et~al.(2014)Wiesemann, Kuhn, \protect\BIBand{}
  Sim}]{wiesemann2014distributionally}
Wiesemann W, Kuhn D, Sim M (2014) Distributionally robust convex optimization.
  \emph{Operations Research} 62(6):1358--1376.

\end{thebibliography}

\newpage

\begin{APPENDICES} \label{sec:appendix}
\section{Technical Proofs}\label{sec:proofs}
Theorem \ref{EL basic} is a simple consequence of the following proposition:
\begin{proposition}
Under the same conditions as Theorem \ref{EL basic}, $\underline Z_n(x)\leq Z_0(x)\leq\overline Z_n(x)$ if and only if $-2\log R(Z_0(x))\leq\chi^2_{1,1-\alpha}$.\label{duality}
\end{proposition}

\proof{Proof of Proposition \ref{duality}.}
We first argue that the optimization defining \eqref{profile1} must have an optimal solution, if it is feasible. Since $-2\sum_{i=1}^n\log(nw_i)\to\infty$ as $w_i\to0$ for any $i$, it suffices to consider only $w_i$ such that $w_i\geq\epsilon$ for some small $\epsilon>0$. Since the set $\left\{\sum_{i=1}^nh(x;\xi_i)w_i=Z_0(x),\ \sum_{i=1}^nw_i=1,\ w_i\geq\epsilon\text{\ for all\ }i=1,\ldots,n\right\}$ is compact, by Weierstrass Theorem, there exists an optimal solution for \eqref{profile1}.

Suppose $-2\log R(Z_0(x))\leq\chi^2_{1,1-\alpha}$. Then the optimization in $-2\log R(Z_0(x))$ is feasible, and there must exist a probability vector $\mathbf w=(w_1,\ldots,w_n)$ such that $-2\sum_{i=1}^n\log(nw_i)\leq\chi^2_{1,1-\alpha}$ and $\sum_{i=1}^nh(x;\xi_i)w_i=Z_0(x)$. This implies $\underline Z_n(x)\leq Z_0(x)\leq\overline Z_n(x)$.

To show the reverse direction, note first that the set
$$\left\{\sum_{i=1}^nh(x;\xi_i)w_i:-2\sum_{i=1}^n\log(nw_i)\leq\chi^2_{1,1-\alpha},\ \sum_{i=1}^nw_i=1,\ w_i\geq0\text{\ for all\ }i=1,\ldots,n\right\}$$
is an interval, since $\sum_{i=1}^nh(x;\xi_i)w_i$ is a linear function of the convex set
$$\left\{(w_1,\ldots,w_n):-2\sum_{i=1}^n\log(nw_i)\leq\chi^2_{1,1-\alpha},\ \sum_{i=1}^nw_i=1,\ w_i\geq0\text{\ for all\ }i=1,\ldots,n\right\}$$
Moreover, since the latter set is compact, by Weierstrass Theorem again, there must exist optimal solutions in the optimization pair
$$\max/\min\left\{\sum_{i=1}^nh(x;\xi_i)w_i:-2\sum_{i=1}^n\log(nw_i)\leq\chi^2_{1,1-\alpha},\ \sum_{i=1}^nw_i=1,\ w_i\geq0\text{\ for all\ }i=1,\ldots,n\right\}$$
Therefore, $\underline Z_n(x)\leq Z_0(x)\leq\overline Z_n(x)$ implies that there exists a probability vector $\mathbf w$ such that $\sum_{i=1}^nh(x;\xi_i)w_i=Z_0(x)$ and $-2\sum_{i=1}^n\log(nw_i)\leq\chi^2_{1,1-\alpha}$, leading to $-2\log R(Z_0(x))\leq\chi^2_{1,1-\alpha}$.\Halmos
\endproof

\proof{Proof of Theorem \ref{EL basic}.}
By Theorem \ref{ELT}, we have $\lim_{n\to\infty}P(-2\log R(Z_0(x))\leq\chi^2_{1,1-\alpha})=1-\alpha$ for a fixed $x\in\Theta$, where
\begin{eqnarray}
&&-2\log R(Z_0(x))\notag\\
&=&\min\left\{-2\sum_{i=1}^n\log(nw_i):\sum_{i=1}^nh(x;\xi_i)w_i=Z_0(x),\ \sum_{i=1}^nw_i=1,\ w_i\geq0\text{\ for all\ }i=1,\ldots,n\right\}\label{profile1}
\end{eqnarray}
Thus, to show \eqref{exact asymptotic}, it suffices to prove that $\underline Z_n(x)\leq Z_0(x)\leq\overline Z_n(x)$ if and only if $-2\log R(Z_0(x))\leq\chi^2_{1,1-\alpha}$. Proposition \ref{duality} finishes the proof.\Halmos
\endproof

\proof{Proof of Proposition \ref{EL discrete}.}
By relabeling the weights under membership of the support points, we rewrite
\begin{eqnarray}
\underline Z_n(x)&=&\min\Bigg\{\sum_{i=1}^kh(x;s_i)\sum_{j=1}^{n_i}w_{ij}:-\frac{1}{n}\sum_{i=1}^k\sum_{j=1}^{n_i}\log(nw_{ij})\leq\frac{\chi^2_{1,1-\alpha}}{2n},{}\notag\\
&&{}\ \sum_{i=1}^k\sum_{j=1}^{n_i}w_{ij}=1,\ w_{ij}\geq0\text{\ for\ }i=1,\ldots,k,\ j=1,\ldots,n_i\Bigg\}\label{opt EL dual min2}
\end{eqnarray}
and
\begin{eqnarray}
\overline Z_n(x)&=&\max\Bigg\{\sum_{i=1}^kh(x;s_i)\sum_{j=1}^{n_i}w_{ij}:-\frac{1}{n}\sum_{i=1}^k\sum_{j=1}^{n_i}\log(nw_{ij})\leq\frac{\chi^2_{1,1-\alpha}}{2n},{}\notag\\
&&{}\ \sum_{i=1}^k\sum_{j=1}^{n_i}w_{ij}=1,\ w_{ij}\geq0\text{\ for\ }i=1,\ldots,k,\ j=1,\ldots,n_i\Bigg\}\label{opt EL dual max2}
\end{eqnarray}
%
%\begin{equation}
%\max\left\{\sum_{j=1}^dh(s_j)\sum_{i=1}^{n_j}w_{ji}:-\frac{1}{n}\sum_{j=1}^d\sum_{i=1}^{n_j}\log(nw_{ji})\leq\frac{\chi^2_{1,1-\alpha}}{2n},\ \sum_{j=1}^d\sum_{i=1}^{n_j}w_{ji}=1,\ w_{ji}\geq0\text{\ for\ }j=1,\ldots,d,\ i=1,\ldots,n_j\right\}\label{opt EL dual max2}
%\end{equation}
%where $w_{ji}$ represents the weight on each data point that falls into $s_j$.
To avoid repetition, we focus on the maximization formulation. We show that, for any feasible $\mathbf p$ in $\max_{\mathbf p\in\mathcal U_{Burg}'}E_{\mathbf p}[h(x;\xi)]$, we can construct a feasible $\mathbf w$ for $\overline Z_n(x)$ that attains the same objective value, and vice versa.

To this end, for any $\mathbf p=(p_1,\ldots,p_k)\in\mathcal U_{Burg}'$, we define $w_{ij}=p_i/n_i$ for all $j=1,\ldots,n_i$. Then $$-\frac{1}{n}\sum_{i=1}^k\sum_{j=1}^{n_i}\log(nw_{ij})=-\sum_{i=1}^k\frac{n_i}{n}\log\frac{np_i}{n_i}=-\sum_{i=1}^k\hat p_i\log\frac{p_i}{\hat p_i}\leq\frac{\chi^2_{1,1-\alpha}}{2n}$$
as well as $\sum_{i=1}^k\sum_{j=1}^{n_i}w_{ij}=\sum_{i=1}^kp_i=1$, and $w_{ij}\geq0$ for all $i$ and $j$. Hence $w_{ij}$ is feasible for \eqref{opt EL dual max2}. Moreover, $\sum_{i=1}^kh(x;s_i)\sum_{j=1}^{n_i}w_{ij}=\sum_{i=1}^kh(x;s_i)p_i$, thus the same objective value is attained.

On the other hand, suppose $\mathbf w=(w_{ij})$ is a feasible solution for \eqref{opt EL dual max2}. We then define $p_i=\sum_{j=1}^{n_i}w_{ij}$. By Jensen's inequality we have $-\log(p_i/n_i)\leq-(1/n_i)\sum_{j=1}^{n_i}\log w_{ij}$, and so
$$-\sum_{i=1}^k\frac{n_i}{n}\log\frac{np_i}{n_i}\leq-\frac{1}{n}\sum_{i=1}^k\sum_{j=1}^{n_i}\log(nw_{ij})\leq\frac{\chi^2_{1,1-\alpha}}{2n}$$
Together with the simple observation that $\sum_{i=1}^kp_i=\sum_{i=1}^k\sum_{j=1}^{n_i}w_{ij}=1$ and $p_i\geq0$ for all $i$, we get that $\mathbf p=(p_i)$ is feasible for $\max_{\mathbf p\in\mathcal U_{Burg}'}E_{\mathbf p}[h(x;\xi)]$. Moreover, $\sum_{i=1}^kh(x;s_i)p_i=\sum_{i=1}^kh(x;s_i)\sum_{j=1}^{n_i}w_{ij}$, thus the same objective value is attained in this case as well.

Similar arguments apply to the minimization formulation, and we conclude the proof. %Therefore \eqref{EL min discrete} and \eqref{EL max discrete} are equivalent to \eqref{opt EL dual min2} and \eqref{opt EL dual max2}.
\Halmos
\endproof

\proof{Proof of Theorem \ref{EL process thm}.}
First, Assumption \ref{finite mean} allows us to define $\tilde h(x;\xi)=h(x;\xi)-Z_0(x)$. Also, we denote the classes of functions $\Xi\to\mathbb R$
$$\mathcal H_\Theta^1=\{|\tilde h(x;\cdot)|:x\in\Theta\}$$
$$\mathcal H_\Theta^2=\{\tilde h(x;\cdot)^2:x\in\Theta\}$$
$$\mathcal H_\Theta^+=\{\tilde h(x;\cdot)^+:x\in\Theta\}$$
$$\mathcal H_\Theta^-=\{\tilde h(x;\cdot)^-:x\in\Theta\}$$
where
$$y^+=\left\{\begin{array}{ll}y&\text{\ \ if\ }y\geq0\\0&\text{\ \ if\ }y<0\end{array}\right.\text{\ \ \ \ and\ \ \ \ }y^-=\left\{\begin{array}{ll}0&\text{\ \ if\ }y>0\\-y&\text{\ \ if\ }y\leq0\end{array}\right.$$
Since $\mathcal H_\Theta$ is a $P_0$-Donsker class, it is $P_0$-Glivenko-Cantelli (GC) (e.g., the discussion before Example 2.1.3 in \cite{van1996weak}). By the preservation theorem (Theorem \ref{preservation} in Appendix \ref{sec:EP}), since $E_0\|\tilde h(\cdot;\xi)^2\|_\Theta=E\|\tilde h(\cdot;\xi)\|_\Theta^2<\infty$ by Assumption \ref{L2}, $\mathcal H_\Theta^2$ is also $P_0$-GC. Moreover, since $E\|\tilde h(\cdot;\xi)^\pm\|_\Theta\leq E\|\tilde h(\cdot;\xi)\|_\Theta\leq\sqrt{E\|\tilde h(\cdot;\xi)\|_\Theta^2}<\infty$, $\mathcal H_\Theta^+$, $\mathcal H_\Theta^-$ and $\mathcal H_\Theta^1$ are all $P_0$-GC as well. Letting $P_n$ be the empirical measure generated from $\xi_1,\ldots,\xi_n$, the above imply
\begin{eqnarray}
\|P_n-P_0\|_{\mathcal H_\Theta^+}\stackrel{a.s.}{\to}0\label{GC+}\\
\|P_n-P_0\|_{\mathcal H_\Theta^-}\stackrel{a.s.}{\to}0\label{GC-}\\
\|P_n-P_0\|_{\mathcal H_\Theta^2}\stackrel{a.s.}{\to}0\label{GC var}
\end{eqnarray}
where $\|P_n-P_0\|_{\mathcal F}=\sup_{f\in\mathcal F}|P_n(f)-P_0(f)|$ for $P_n$ indexed by $\mathcal F$, and similarly defined for $P_0$ (see Appendix \ref{sec:EP}).

%\sup_{\theta\in\Theta}|\mathbb P_n(\tilde h(x;\cdot)^2)-P(\tilde h(x;\cdot)^2)|\stackrel{a.s.}{\to}0\label{GC var}
%\sup_{\theta\in\Theta}|\mathbb P_n(\tilde h(x;\cdot)^+)-P(\tilde h(x;\cdot)^+)|\stackrel{a.s.}{\to}0\label{GC+}\\
%\sup_{\theta\in\Theta}|\mathbb P_n(\tilde h(x;\cdot)^-)-P(\tilde h(x;\cdot)^-)|\stackrel{a.s.}{\to}0\label{GC-}
 %a $P$ class $
%$$\{(h(\cdot;\theta)-Z(\theta))^2:\theta\in\Theta\}\text{\ \ and\ \ }\{(h(\cdot;\theta)-Z(\theta))^-:\theta\in\Theta\}$$
%and
%$$\sup_{\theta\in\Theta}|\mathcal P_n((h(\cdot;\theta)-Z(\theta))^-)-P((h(\cdot;\theta)-Z(\theta))^-)|\stackrel{a.s.}{\to}0$$
%or equivalently
%\begin{equation}
%\sup_{\theta\in\Theta}|\mathbb V_n(h(x))-Var(h(x;\xi))|\stackrel{a.s.}{\to}0\label{convergence var}
%\end{equation}
%where $\mathbb V_n(h(x))=\mathbb P_n((h(x;\cdot)-Z(x))^2)$ is the empirical variance (with known mean function $Z(x)$).
Note that \eqref{GC var} in particular implies the uniform convergence of the empirical variance
\begin{equation}
\left\|\frac{1}{n}\sum_{i=1}^n(h(\cdot;\xi_i)-Z_0(\cdot))^2-\sigma_0(\cdot)\right\|_\Theta\stackrel{a.s.}{\to}0\label{convergence var}
\end{equation}
where $\sigma_0(x)=Var_0(h(x;\xi))$.

Now, for each $x$,
$$E_0\tilde h(x;\xi)^++E_0\tilde h(x;\xi)^-=E_0|\tilde h(x;\xi)|$$
$$E_0\tilde h(x;\xi)^+-E_0\tilde h(x;\xi)^-=E_0\tilde h(x;\xi)=0$$
which gives $E_0\tilde h(x;\xi)^+=E_0\tilde h(x;\xi)^-=E_0|\tilde h(x;\xi)|/2$. Hence
\begin{equation}
\inf_{x\in\Theta}E_0\tilde h(x;\xi)^\pm=\inf_{x\in\Theta}\frac{E_0|\tilde h(x;\xi)|}{2}\geq\frac{c}{2}\label{lower bdd pm}
\end{equation}
where we define $c$ as a constant such that $\inf_{x\in\Theta}E_0|h(x;\xi)-Z_0(x)|\geq c$, which exists by Assumption \ref{nondegeneracy}. By Jensen's inequality,
\begin{equation}
\inf_{x\in\Theta}E_0\tilde h(x;\xi)^2\geq\left(\inf_{x\in\Theta}E_0|\tilde h(x;\xi)|\right)^2\geq c^2\label{lower bdd var}
\end{equation}
by using Assumption \ref{nondegeneracy} again. From \eqref{GC+}, \eqref{GC-} and \eqref{lower bdd pm}, we have $\inf_{x\in\Theta}(1/n)\sum_{i=1}^n\tilde h(x;\xi_i)^+$ and $\inf_{x\in\Theta}(1/n)\sum_{i=1}^n\tilde h(x;\xi_i)^->0$ for large enough $n$ a.s.. When this occurs, $\min_{1\leq i\leq n}h(x;\xi_i)<Z_0(x)<\max_{1\leq i\leq n}h(x;\xi_i)$ for every $x$, and the optimization defining $-\log R(x;Z_0)$, namely
\begin{equation}
\min\left\{-\sum_{i=1}^n\log(nw_i)\Bigg|\sum_{i=1}^nh(x;\xi_i)w_i=Z_0(x),\ \sum_{i=1}^nw_i=1,\ w_i\geq0\text{\ for\ }i=1,\ldots,n\right\}\label{EL opt transformed}
\end{equation}
has a unique optimal solution $\mathbf w(x)=(w_1(x),\ldots,w_n(x))$ with $w_i(x)>0$ for all $i$, for any $x$. This is because setting any $w_i(x)=0$ would render $-2\sum_{i=1}^n\log(nw_i)=\infty$ which is clearly suboptimal. Hence it suffices to replace $w_i\geq0$ with $w_i\geq\epsilon$ for all $i$ for some small enough $\epsilon>0$. In this modified region, the optimum exists and is unique since $-\sum_{i=1}^n\log(nw_i)$ is strictly convex.%Moreover, this solution satisfies the KKT condition (REF).

%Since $\inf_{\theta\in\Theta}Var(h(X;\theta))\geq c>0$ by the assumption, we have $\inf_{\theta\in\Theta}\mathbb V_n(h(X;\theta))>0$ as for large enough $n$ a.s..

Now consider the optimization \eqref{EL opt transformed} when $\min_{1\leq i\leq n}h(x;\xi_i)<Z_0(x)<\max_{1\leq i\leq n}h(x;\xi_i)$. We adopt the proof technique in Section 11.2 of \cite{owen2001empirical}, but generalize at the empirical process level.
For convenience, we write $\tilde h_i=\tilde h(x;\xi_i)=h(x;\xi_i)-Z_0(x)$, and also suppress the $x$ in $w_i=w_i(x)$ and $Z_0=Z_0(x)$. The Lagrangian is written as
$$-\sum_{i=1}^n\log(nw_i)+\lambda\left(\sum_{i=1}^n\tilde h_iw_i-Z_0\right)+\gamma\left(\sum_{i=1}^nw_i-1\right)$$
where $\lambda=\lambda(x)$ and $\gamma=\gamma(x)$ are the Lagrange multipliers.
Differentiating with respect to $w_i$ and setting it to zero, we have
\begin{equation}
-\frac{1}{w_i}+\lambda\tilde h_i+\gamma=0\label{interim}
\end{equation}
Setting $\sum_{i=1}^n\tilde h_iw_i=0$ and $\sum_{i=1}^nw_i=1$, multiplying both sides of \eqref{interim} by $w_i$ and summing up over $i$, we get $\gamma=n$. Using \eqref{interim} again, we have
\begin{equation}
w_i=\frac{1}{n}\frac{1}{1+\lambda\tilde h_i}\label{w}
\end{equation}
where the $\lambda$ in \eqref{w} is rescaled by a factor of $n$. Note that we can find $\lambda$ such that
\begin{equation}
\sum_{i=1}^n\frac{1}{n}\frac{\tilde h_i}{1+\lambda\tilde h_i}=0\label{interim3}
\end{equation}
and $\frac{1}{n}\frac{1}{1+\lambda\tilde h_i}>0$ for all $i$, upon which the Karush-Kuhn-Tucker (KKT) condition can be seen to hold and conclude that $w_i$ in \eqref{w} is the optimal solution. Indeed, let $\tilde h^*=\max_i\tilde h_i>0$ and $\tilde h_*=\min_i\tilde h_i<0$. Note that $\sum_{i=1}^n\frac{1}{n}\frac{\tilde h_i}{1+\lambda\tilde h_i}\to\infty$ as $\lambda\to-1/\tilde h^*$, and $\to-\infty$ as $\lambda\to-1/\tilde h_*$. Since $\sum_{i=1}^n\frac{1}{n}\frac{\tilde h_i}{1+\lambda\tilde h_i}$ is a continuous function in $\lambda$ between $-1/\tilde h^*$ and $-1/\tilde h_*$, there must exist a $\lambda$ that solves \eqref{interim3}. Moreover, $\frac{1}{n}\frac{1}{1+\lambda\tilde h_i}>0$ for all $i$ for this $\lambda$.

Given this characterization of the optimal solution, the rest of the proof is to derive the asymptotic behavior of $-2\log R(x;Z_0)$ as $n\to\infty$. First, we write
$$\frac{1}{1+\lambda\tilde h_i}=1-\frac{\lambda\tilde h_i}{1+\lambda\tilde h_i}$$
Multiplying both sides by $\tilde h_i/n$ and summing up over $i$, we get
$$\sum_{i=1}^n\frac{1}{n}\frac{\tilde h_i}{1+\lambda\tilde h_i}=\bar h-\lambda\sum_{i=1}^n\frac{1}{n}\frac{\tilde h_i^2}{1+\lambda\tilde h_i}$$
where $\bar h:=\bar h(x)=(1/n)\sum_{i=1}^n\tilde h(x;\xi_i)$, and hence
\begin{equation}
\bar h=\lambda\sum_{i=1}^n\frac{1}{n}\frac{\tilde h_i^2}{1+\lambda\tilde h_i}\label{interim2}
\end{equation}
by \eqref{interim3}. Now let
$$s:=s(x)=\frac{1}{n}\sum_{i=1}^n\tilde h(x;\xi_i)^2$$
be the empirical variance. Note that, from \eqref{w}, $w_i>0$ implies $1+\lambda\tilde h_i>0$. Together with $s\geq0$, we get
\begin{align*}
|\lambda|s&\leq\left|\lambda\sum_{i=1}^n\frac{1}{n}\frac{\tilde h_i^2}{1+\lambda\tilde h_i}\right|\left(1+|\lambda|\max_{1\leq i\leq n}|\tilde h_i|\right)\\
&=|\bar h|\left(1+|\lambda|\max_{1\leq i\leq n}|\tilde h_i|\right)
\end{align*}
by using \eqref{interim2}, and hence
\begin{equation}
|\lambda|\left(s-|\bar h|\max_{1\leq i\leq n}|\tilde h_i|\right)\leq|\bar h|\label{interim new}
\end{equation}
By Lemma \ref{lemma max} in Appendix \ref{sec:thms}, $E\|\tilde h(x;\xi)\|_\Theta^2<\infty$ in Assumption \ref{L2} implies that
\begin{equation}
\max_{1\leq i\leq n}\|\tilde h_i\|_\Theta=o(n^{1/2})\text{\ a.s.}\label{interim new1}
\end{equation}
Moreover, since $\mathcal H_\Theta$ is $P_0$-Donsker, we have $\sqrt n\bar h=\sqrt n(P_n(h(\cdot;\cdot))-P_0(h(\cdot;\cdot)))\Rightarrow\tilde G$ in $\ell^\infty(\mathcal H_\Theta)$, where $\tilde G(\cdot)$ is a tight Gaussian process indexed by $h(x;\cdot)\in\mathcal H_\Theta$ that is centered and has covariance $Cov(\tilde G(h(x_1;\cdot)),\tilde G(h(x_2;\cdot)))=Cov_0(h(x_1;\xi),h(x_2;\xi))$ for any $h(x_1;\cdot),h(x_2;\cdot)\in\mathcal H_\Theta$. Noting that the map $\ell^\infty(\mathcal H_\Theta)\to\ell^\infty(\Theta)$ defined by $y(\cdot)\mapsto y(h(\cdot;\cdot))$ is continuous, by the continuous mapping theorem (Theorem \ref{continuous mapping} in Appendix \ref{sec:thms}), we have $\sqrt n\bar h\Rightarrow\tilde G$ in $\ell^\infty(\Theta)$ where $\tilde G$ is now indexed by $x\in\Theta$. As the norm in $\ell^\infty(\Theta)$, $\|\cdot\|_\Theta$ is a continuous map. By the continuous mapping theorem again, $\sqrt n\|\bar h\|_\Theta=\|\sqrt n\bar h\|_\Theta\Rightarrow\|\tilde G\|_\Theta$, %Since $\|\bar h\|_{\mathcal H_\Theta}=\|\bar h\|_\Theta$ and $\|\tilde G\|_{\mathcal H_\Theta}=\|\tilde G\|_\Theta$, we have $\sqrt\|\bar h\|_\Theta\Rightarrow\|\tilde G\|_\Theta$ and
so that $\|\bar h\|_\Theta=O_p(n^{-1/2})$. Moreover, $\||\bar h|\max_{1\leq i\leq n}|\tilde h_i|\|_\Theta\leq\|\bar h\|_\Theta\max_{1\leq i\leq n}\|\tilde h_i\|_\Theta=O_p(n^{-1/2})o(n^{1/2})=o_p(1)$.

Next, from \eqref{convergence var} and \eqref{lower bdd var}, we have $\inf_{x\in\Theta}s(x)\geq c_1$ for some $c_1>0$ ev.. Pick any constant $\varepsilon<c_1$. We have
\begin{equation}
P\left(\inf_{x\in\Theta}\left\{s(x)-|\bar h(x)|\max_{1\leq i\leq n}|\tilde h_i(x)|\right\}\geq c_1-\varepsilon\right)\geq P\left(\inf_{x\in\Theta}s(x)\geq c_1,\left\||\bar h(x)|\max_{1\leq i\leq n}|\bar h_i(x)|\right\|_\Theta\leq\varepsilon\right)\to1\label{interim11}
\end{equation}
Over the set $\{\inf_{x\in\Theta}\{s(x)-|\bar h(x)|\max_{1\leq i\leq n}|\bar h_i(x)|\}>c_1-\varepsilon\}$, \eqref{interim new} implies
$$|\lambda(x)|\leq\frac{|\bar h(x)|}{s-|\bar h(x)|\max_{1\leq i\leq n}|\bar h_i(x)|}$$
for all $x\in\Theta$, so that
\begin{equation}
\|\lambda\|_\Theta\leq\frac{\|\bar h\|_\Theta}{c_1-\varepsilon}\label{interim updated1}
\end{equation}
We argue that $\|\lambda\|_\Theta=O_p(n^{-1/2})$. This is because, for any given $\delta>0$, we can find a large enough $B>0$ such that
\begin{eqnarray*}
&&\limsup_{n\to\infty}P(\|\lambda\|_\Theta>Bn^{-1/2})\\
&\leq&\limsup_{n\to\infty}\Bigg(P\left(\|\lambda\|_\Theta>Bn^{-1/2},\ \inf_{x\in\Theta}\left\{s(x)-|\bar h(x)|\max_{1\leq i\leq n}|\tilde h_i(x)|\right\}>c_1-\varepsilon\right){}\\
&&{}+P\left(\inf_{x\in\Theta}\left\{s(x)-|\bar h(x)|\max_{1\leq i\leq n}|\tilde h_i(x)|\right\}\leq c_1-\varepsilon\right)\Bigg)\\
&\leq&\limsup_{n\to\infty}\left(P(\|\bar h\|_\Theta>(c_1-\varepsilon)Bn^{-1/2})+P\left(\inf_{x\in\Theta}\left\{s(x)-|\bar h(x)|\max_{1\leq i\leq n}|\tilde h_i(x)|\right\}\leq c_1-\varepsilon\right)\right){}\\
&&{}\text{\ \ by \eqref{interim updated1}}\\
&<&\delta
\end{eqnarray*}
by \eqref{interim11} and that $\|\bar h\|_\Theta=O_p(n^{-1/2})$ as shown above. This and \eqref{interim new1} together gives
\begin{equation}
\max_{1\leq i\leq n}\sup_{x\in\Theta}|\lambda(x)\tilde h_i(x)|\leq\|\lambda\|_\Theta\max_{1\leq i\leq n}\|\tilde h_i\|_\Theta=O_p(n^{-1/2})o(n^{1/2})=o_p(1)\label{interim6}
\end{equation}

Now \eqref{interim3} can be rewritten as
\begin{align}
0&=\sum_{i=1}^n\frac{1}{n}\tilde h_i\left(1-\lambda\tilde h_i+\frac{\lambda^2\tilde h_i^2}{1+\lambda\tilde h_i}\right)\notag\\
&=\bar h-\lambda s+\sum_{i=1}^n\frac{1}{n}\frac{\lambda^2\tilde h_i^3}{1+\lambda\tilde h_i}\label{interim4}
\end{align}
The last term in \eqref{interim4} satisfies
$$\left|\sum_{i=1}^n\frac{1}{n}\frac{\lambda^2\tilde h_i^3}{1+\lambda\tilde h_i}\right|\leq\frac{1}{n}\sum_{i=1}^n\tilde h_i^2\lambda^2\max_{1\leq i\leq n}|\tilde h_i|\max_{1\leq i\leq n}(1+\lambda\tilde h_i)^{-1}$$
Taking $\sup_{\theta\in\Theta}$ on both sides, we get
\begin{equation}
\left\|\sum_{i=1}^n\frac{1}{n}\frac{\lambda^2\tilde h_i^3}{1+\lambda\tilde h_i}\right\|_\Theta\leq\|s\|_\Theta\|\lambda\|_\Theta^2\max_{1\leq i\leq n}\|\tilde h_i\|_\Theta\max_{1\leq i\leq n}\sup_{x\in\Theta}(1+\lambda\tilde h(x,\xi_i))^{-1}\label{interim new2}
\end{equation}
Now, $\|s\|_\Theta\to\|\sigma_0\|_\Theta$ by \eqref{convergence var}. Moreover, for any small $\varepsilon>0$,
\begin{align*}
P\left(\max_{1\leq i\leq n}\sup_{x\in\Theta}(1+\lambda\tilde h(x,\xi_i))^{-1}>\frac{1}{1-\varepsilon}\right)&=P\left(\frac{1}{1+\lambda\tilde h(x,\xi_i)}>\frac{1}{1-\varepsilon}\text{\ \ for some $1\leq i\leq n$ and $\theta\in\Theta$}\right)\\
&\leq P(\lambda\tilde h(x,\xi_i)<-\varepsilon\text{\ \ for some $1\leq i\leq n$ and $x\in\Theta$})\\
&\leq P\left(\max_{1\leq i\leq n}\sup_{x\in\Theta}|\lambda\tilde h(x,\xi_i)|>\varepsilon\right)\to0
\end{align*}
by \eqref{interim6}. Thus \eqref{interim new2} is bounded from above by
\begin{equation}
O(1)O_p(n^{-1})o(n^{1/2})O_p(1)=o_p(n^{-1/2})\label{interim5}
\end{equation}

From \eqref{interim4} and \eqref{interim5}, we have
$$0=\bar h-\lambda s+\epsilon$$
where $\|\epsilon\|_\Theta=o_p(n^{-1/2})$. Since \eqref{convergence var} and \eqref{lower bdd var} implies $\inf_{x\in\Theta}s(x)\geq c_1$ for some $c_1>0$ ev., we further get
\begin{equation}
\lambda=s^{-1}(\bar h+\epsilon)\label{interim new4}
\end{equation}
%where here $\epsilon$ is another error term that again satisfies $\|\epsilon\|_\Theta=o_p(n^{-1/2})$.

Now consider
\begin{align}
-2\log R(x;Z_0)&=-2\sum_{i=1}^n\log(nw_i)\notag\\
&=2\sum_{i=1}^n\log(1+\lambda\tilde h_i)\notag\\
&=2\sum_{i=1}^n\left(\lambda\tilde h_i-\frac{1}{2}\lambda^2\tilde h_i^2+\nu_i\right)\label{interim7}
\end{align}
where
$$\nu_i=\frac{1}{3}\frac{1}{(1+\zeta_i)^3}(\lambda\tilde h_i)^3$$
with $\zeta_i=\zeta_i(x)$ between 0 and $\lambda\tilde h_i$ by Taylor's expansion. So
\begin{equation}
|\nu_i|\leq\frac{1}{3}\frac{1}{|1+\zeta_i|^3}|\lambda\tilde h_i|^3\label{interim new3}
\end{equation}
For any large enough $B_1>0$, we have
\begin{eqnarray}
&&P(|\nu_i(x)|>B_1|\lambda(x)\tilde h(x,\xi_i)|^3\text{\ \ for all $x\in\Theta$ and $1\leq i\leq n$})\notag\\
&\leq&P\left((3B_1)^{1/3}\max_{1\leq i\leq n}\|1+\zeta_i\|_\Theta<1\right)\text{\ \ from \eqref{interim new3}}\notag\\
&=&P\left((3B_1)^{1/3}\max_{1\leq i\leq n}\|1+\zeta_i\|_\Theta<1,\max_{1\leq i\leq n}\|\zeta_i\|_\Theta<\varepsilon\right){}\notag\\
&&{}+P\left((3B_1)^{1/3}\max_{1\leq i\leq n}\|1+\zeta_i\|_\Theta<1,\max_{1\leq i\leq n}\|\zeta_i\|_\Theta>\varepsilon\right)\text{\ \ for some sufficiently large $0<\varepsilon<1$}\notag\\
&\leq&P\left((3B_1)^{1/3}\max_{1\leq i\leq n}(1-\varepsilon)<1\right)+P\left(\max_{1\leq i\leq n}\|\zeta_i\|_\Theta>\varepsilon\right)\notag\\
&\to&0\label{interim13}
\end{eqnarray}
since $\max_{1\leq i\leq n}\|\zeta_i\|_\Theta\leq\max_{1\leq i\leq n}\|\lambda\tilde h_i\|_\Theta=o_p(1)$ by \eqref{interim6}.
%\begin{equation}
%P(|\nu_i|\leq B_1|\lambda\tilde h_i(x)|^3\text{\ \ for all $x\in\Theta$ and $1\leq i\leq n$})\to1\label{interim13}
%\end{equation}
%for some $B_1>0$.
Now \eqref{interim7} gives
\begin{eqnarray}
&&2n\lambda\bar h-\lambda^2ns+2\sum_{i=1}^n\nu_i\notag\\
&=&2ns^{-1}(\bar h+\epsilon)\bar h-nss^{-2}(\bar h^2+2\epsilon\bar h+\epsilon^2)+2\sum_{i=1}^n\nu_i\text{\ \ by \eqref{interim new4}}\notag\\
&=&ns^{-1}\bar h^2-ns^{-1}\epsilon^2+2\sum_{i=1}^n\nu_i\label{interim10}
\end{eqnarray}
Note that
\begin{equation}
\|ns^{-1}\epsilon^2\|_\Theta\leq nO(1)o_p(n^{-1})=o_p(1)\label{interim8}
\end{equation}
since $\inf_{x\in\Theta}s(x)\geq c_1>0$ ev., and $\|\epsilon\|_\Theta=o_p(n^{-1/2})$. Moreover, over the set $\{|\nu_i(x)|\leq B_1|\lambda(x)\tilde h_i(x)|^3$ for all $x\in\Theta$ and $1\leq i\leq n\}$,
we have
$$\left|\sum_{i=1}^n\nu_i\right|\leq B_1|\lambda|^3\sum_{i=1}^n|\tilde h_i|^2\max_{1\leq i\leq n}|\tilde h_i|$$
for all $x\in\Theta$. Note that
\begin{equation}
\left\|B_1|\lambda|^3\sum_{i=1}^n|\tilde h_i|^2\max_{1\leq i\leq n}|\tilde h_i|\right\|_\Theta\leq B_1O_p(n^{-3/2})nO(1)o(n^{1/2})=o_p(1)\label{interim12}
\end{equation}
since $\|\lambda\|_\Theta=O_p(n^{-1/2})$, $\|s\|_\Theta=O(1)$, and $\max_{1\leq i\leq n}\|\tilde h_i\|_\Theta=o(n^{1/2})$. Now, for any $\varepsilon>0$,
\begin{eqnarray*}
&&P\left(\left\|\sum_{i=1}^n\nu_i\right\|_\Theta>\varepsilon\right)\\
&\leq&P\left(\left\|\sum_{i=1}^n\nu_i\right\|_\Theta>\varepsilon,\ |\nu_i(x)|\leq B_1|\lambda(x)\tilde h_i(x)|^3\text{\ \ for some $x\in\Theta$ and $1\leq i\leq n$}\right){}\\
&&{}+P(|\nu_i(x)|>B_1|\lambda(x)\tilde h_i(x)|^3\text{\ \ for all $x\in\Theta$ and $1\leq i\leq n$})\\
&\leq&P\left(\left\|B_1|\lambda|^3\sum_{i=1}^n|\tilde h_i|^2\max_{1\leq i\leq n}|\tilde h_i|\right\|_\Theta>\varepsilon\right)+P(|\nu_i|>B_1|\lambda\tilde h_i|^3\text{\ \ for all $x\in\Theta$ and $1\leq i\leq n$})\\
&\to&0
\end{eqnarray*}
by \eqref{interim13} and \eqref{interim12}. Hence we have
\begin{equation}
\left\|\sum_{i=1}^n\nu_i\right\|_\Theta=o_p(1)\label{interim9}
\end{equation}
Using \eqref{interim8} and \eqref{interim9}, \eqref{interim10} implies that
\begin{equation}
-2\log R(x;Z_0)=ns^{-1}\bar h^2+\epsilon_1\label{interim new5}
\end{equation}
where $\|\epsilon_1\|_\Theta=o_p(1)$.

Note that $\sqrt n\bar h\Rightarrow\tilde G$ in $\ell^\infty(\Theta)$ where $\tilde G(\cdot)$ is defined previously as the centered Gaussian process with covariance $Cov(\tilde G(x_1),\tilde G(x_2))=Cov(h(x_1;\xi),h(x_2;\xi))$ for any $x_1,x_2\in\Theta$. By Slutsky's Theorem (Theorem \ref{Slutsky} in Appendix \ref{sec:thms}) and \eqref{convergence var}, $(\sqrt n\bar h,s)\Rightarrow(\tilde G,\sigma_0)$ in $(\ell^\infty\times\ell^\infty)(\Theta)$ defined as
$$(\ell^\infty\times\ell^\infty)(\Theta)=\left\{(y_1,y_2):\Theta\to\mathbb R^2\Bigg|\|y_1\|_\Theta+\|y_2\|_\Theta<\infty\right\}$$
Note that pointwise division and $(\cdot)^2$ are continuous maps on $(\ell^\infty\times\ell^\infty)(\Theta)\to\ell^\infty(\Theta)$ and $\ell^\infty(\Theta)\to\ell^\infty(\Theta)$ respectively. Also, $\inf_{x\in\Theta}\sigma_0(x)>0$ by Assumption \ref{nondegeneracy} and Jensen's inequality. By continuous mapping theorem, we have $ns^{-1}\bar h^2\Rightarrow J$ in $\ell^\infty(\Theta)$ where $J(\cdot)$ is as defined in the theorem. % centered and covariance
%$$Cov(G(x_1),G(x_2))=\frac{Cov(h(x_1;\xi),h(x_2;\xi))}{\sqrt{Var(h(x_1;\xi))Var(h(x_2;\xi))}}$$
%for $x_1,x_2\in\Theta$. In particular, $Var(G(x))=1$ for any $x\in\Theta$.
Finally, from \eqref{interim new5} and $\|\epsilon_1\|_\Theta=o_p(1)$, we get further that $-2\log R(\cdot;Z_0)\Rightarrow J(\cdot)$ in $\ell^\infty(\Theta)$. This concludes the proof.\Halmos
\endproof

\proof{Proof of Lemma \ref{sample Gaussian}.}
First define $\tilde G(\cdot)$ as a centered Gaussian process indexed by $x\in\Theta$ with covariance $Cov(\tilde G(x_1),\tilde G(x_2))=Cov_0(h(x_1;\xi),h(x_2;\xi))$. Conditional on almost every data realization $(P_n:n\geq1)$, define $\tilde G_n(\cdot)$ as a centered Gaussian process indexed by $x\in\Theta$ with covariance $Cov(\tilde G_n(x_1),\tilde G_n(x_2))=(1/n)\sum_{i=1}^n(h(x_1;\xi)-\hat h(x_1))(h(x_2;\xi)-\hat h(x_2))$, and $\hat h(x)=(1/n)\sum_{i=1}^nh(x;\xi_i)$.

We first show that $\tilde G_n(\cdot)\Rightarrow\tilde G(\cdot)$ in $\ell^\infty(\Theta)$. Note that, by the property of Gaussian processes, any finite-dimensional vector $(\tilde G_n(x_1),\ldots,\tilde G_n(x_d))$ is distributed as $N(\mathbf 0,\Sigma_n)$, where $\mathbf 0$ is the zero vector and
$$\Sigma_n=\left(\frac{1}{n}\sum_{k=1}^n(h(x_i;\xi_k)-\hat h(x_i))(h(x_j;\xi_k)-\hat h(x_j))\right)_{i,j=1,\ldots,d}$$
%$\Sigma_n^{1/2}$ is the lower triangular square-root matrix of $\Sigma$ obtained from Cholesky's factorization, and $\mathbf Z\sim N(0,I_d)$ for an $d$-dimensional identity matrix $I_d$.

On the other hand, $(\tilde G(x_1),\ldots,\tilde G(x_d))$ is distributed as $N(\mathbf 0,\Sigma)$, where $\Sigma=(Cov_0(h(x_i;\xi),h(x_j;\xi))_{i,j=1,\ldots,d}$. Note that $\Sigma_n\to\Sigma$ a.s. in each entry, and hence $(\tilde G_n(x_1),\ldots,\tilde G_n(x_d))\Rightarrow(\tilde G(x_1),\ldots,\tilde G(x_d))$ (by using for example convergence of the characteristic function). %Moreover, by Assumption \ref{nondegeneracy}, $\Sigma$ is positive definite, and hence $\Sigma^{1/2}$ is unique and $\Sigma_n^{1/2}\to\Sigma$ a.s. in each entry. By Slutsky's Theorem, $\Sigma_n^{1/2}\mathbf Z\Rightarrow\Sigma^{1/2}\mathbf Z$, which implies that any finite-dimensional projection of $\tilde G_n(\cdot)$ converges in distribution to that of $\tilde G(\cdot)$.

Next, note that by Assumption \ref{complexity}, $\mathcal H_\Theta$ is $P_0$-Donsker and hence is totally bounded equipped with the semi-metric $\rho_0(h(x_1;\cdot),h(x_2;\cdot)):=(Var_0(h(x_1;\xi)-h(x_2;\xi)))^{1/2}$ (Section 2.1.2 in \cite{van1996weak}). Equivalently, $\Theta$ is totally bounded under the semi-metric $\rho_0(x_1,x_2):=(Var_0((h(x_1;\xi)-h(x_2;\xi)))^{1/2}$.

We shall also show that $\tilde G_n(\cdot)$ is uniformly equicontinuous in probability under the same semi-metric. To this end, we want to show that
\begin{equation}
\lim_{\delta\to0}\limsup_{n\to\infty}P_{\bm\xi}\left(\sup_{\rho_0(x_1x_2)<\delta}|\tilde G_n(x_1)-\tilde G_n(x_2)|>\epsilon\right)=0\label{asymptotic equicontinuity}
\end{equation}
where $P_{\bm\xi}(\cdot)$ is the probability conditional on the data $\bm\xi$. First, by using the covariance structure of the Gaussian process $\tilde G_n(\cdot)$,
\begin{eqnarray*}
&&E(\tilde G_n(x_1)-\tilde G_n(x_2))^2=\widehat{Var}_n(h(x_1;\xi)-h(x_2;\xi))\\
&=&\frac{1}{n}\sum_{i=1}^n(h(x_1;\xi)-\hat h(x_1))^2+\frac{1}{n}\sum_{i=1}^n(h(x_2;\xi)-\hat h(x_2))^2-\frac{2}{n}\sum_{i=1}^n(h(x_1;\xi)-\hat h(x_1))(h(x_2;\xi)-\hat h(x_2))
\end{eqnarray*}

%
%
%\notag\\
%&=&\lim_{\delta\to0}\limsup_{n\to\infty}P_{\bm\xi}\left(\sup_{\rho_0(x_1,x_2)<\delta}\left(\widehat{Var}_n(h(x_1;\xi)-h(x_2;\xi))\right)^{1/2}|Z|>\epsilon\right)\notag\\
%&=&\lim_{\delta\to0}\limsup_{n\to\infty}2\bar\Phi\left(\frac{\epsilon}{\sup_{\rho_0(x_1,x_2)<\delta}\left(\widehat{Var}_n(h(x_1;\xi)-h(x_2;\xi))\right)^{1/2}}\right)\label{interim new21}
%\end{eqnarray}
%where $Z\sim N(0,1)$, $\bar\Phi(x)$ is the tail distribution of $N(0,1)$, and, by using

Now define $\tilde h(x;\xi)=h(x;\xi)-Z_0(x)$ under Assumption \ref{finite mean}. Note that
\begin{align*}
E\left[\sup_{x,y\in\Theta}(\tilde h(x;\xi)-\tilde h(y;\xi))^2\right]&\leq E\left[\sup_{x\in\Theta}\tilde h(x;\xi)^2+\sup_{y\in\Theta}\tilde h(y;\xi)^2+2\sup_{x\in\Theta}|\tilde h(x;\xi)|\sup_{y\in\Theta}|\tilde h(x;\xi)|\right]\\
&=4E\|\tilde h(\cdot;\xi)\|_\Theta^2<\infty
\end{align*}
by Assumption \ref{L2}. Viewing $\tilde h(x;\cdot)$ and $\tilde h(y;\cdot)$ each as a function $(x,y)\in\Theta^2\to\mathbb R$, we can apply the preservation theorem to conclude that the class of functions
$$\mathcal H_\Theta^\Pi=\{(\tilde h(x;\cdot)-\tilde h(y;\cdot))^2:(x,y)\in\Theta^2\}$$
is a $P_0$-GC class. Therefore,
\begin{equation}
\sup_{x,y\in\Theta}\left|\frac{1}{n}\sum_{i=1}^n(\tilde h(x;\xi_i)-\tilde h(y;\xi_i))^2-E_0(\tilde h(x;\xi_i)-\tilde h(y;\xi_i))^2\right|\to0\text{\ \ a.s.}\label{interim new22}
\end{equation}
Now, note that
\begin{eqnarray}
&&\widehat{Var}_n(h(x_1;\xi)-h(x_2;\xi))\notag\\
&=&\frac{1}{n}\sum_{i=1}^n\left(\left(h(x_1;\xi_i)-\hat h(x_1)\right)-\left(h(x_2;\xi_i)-\hat h(x_2)\right)\right)^2\notag\\
&=&\frac{1}{n}\sum_{i=1}^n\left(\left(\tilde h(x_1;\xi_i)-\tilde h(x_2;\xi_i)\right)-\left((\hat h(x_1)-Z_0(x_1))-(\hat h(x_2)-Z_0(x_2))\right)\right)^2\notag\\
&=&\frac{1}{n}\sum_{i=1}^n\left(\tilde h(x_1;\xi_i)-\tilde h(x_2;\xi_i)\right)^2-\left(\left(\hat h(x_1)-Z_0(x_1)\right)-\left(\hat h(x_2)-Z_0(x_2)\right)\right)^2\label{interim updated2}
\end{eqnarray}
Since $\mathcal H_\Theta$ is $P_0$-GC, $\|\hat h(\cdot)-Z_0(\cdot)\|_\Theta\to0$ a.s.. Hence $\sup_{x_1,x_2\in\Theta}\left(\left(\hat h(x_1)-Z_0(x_1)\right)-\left(\hat h(x_2)-Z_0(x_2)\right)\right)^2\to0$ a.s.. Combining with \eqref{interim new22}, we have, from \eqref{interim updated2},
\begin{equation}
\sup_{x_1,x_2\in\Theta}|\widehat{Var}_n(h(x_1;\xi)-h(x_2;\xi))-Var_0(h(x_1,\xi)-h(x_2,\xi))|\to0\text{\ \ a.s.}\label{interim new23}
\end{equation}
by noting that $Var_0(h(x_1,\xi)-h(x_2,\xi))=E_0(\tilde h(x;\xi)-\tilde h(y;\xi))^2$.

Therefore, by \eqref{interim new23}, for any $\gamma>1$, we have
$$E_{\bm\xi}(\tilde G_n(x_1)-\tilde G_n(x_2))^2=\widehat{Var}_n(h(x_1;\xi)-h(x_2;\xi))\leq\gamma Var_0(h(x_1;\xi)-h(x_2;\xi))=E(\tilde G_\gamma(x_1)-\tilde G_\gamma(x_2))^2$$
a.s. for any $x_1,x_2\in\Theta$, when $n$ is sufficiently large, where $\tilde G_\gamma(\cdot):=\gamma G(\cdot)$ and $E_{\bm\xi}[\cdot]$ denotes the expectation conditional on $\bm\xi$. Thus, by the argument for the Sudakov-Fernique inequality (the first equation in the proof of Theorem 2.9 in \cite{adler1990introduction}), we have
$$E_{\bm\xi}\left[\sup_{\rho_0(x_1,x_2)<\delta}|\tilde G_n(x_1)-\tilde G_n(x_2)|\right]\leq E\left[\sup_{\rho_0(x_1,x_2)<\delta}|\tilde G_\gamma(x_1)-\tilde G_\gamma(x_2)|\right]=\gamma E\left[\sup_{\rho_0(x_1,x_2)<\delta}|\tilde G(x_1)-\tilde G(x_2)|\right]$$
when $n$ is large. Note that
$$\lim_{\delta\to0}E\left[\sup_{\rho_0(x_1,x_2)<\delta}|\tilde G(x_1)-\tilde G(x_2)|\right]=0$$
since $\tilde G(\cdot)$ is tight by the $P_0$-Donsker property of $\mathcal H_\Theta$. Thus
\begin{eqnarray*}
&&\limsup_{n\to\infty}P_{\bm\xi}\left(\sup_{\rho_0(x_1x_2)<\delta}|\tilde G_n(x_1)-\tilde G_n(x_2)|>\epsilon\right)\\
&\leq&\limsup_{n\to\infty}\frac{E_{\bm\xi}\left[\sup_{\rho_0(x_1x_2)<\delta}|\tilde G_n(x_1)-\tilde G_n(x_2)|\right]}{\epsilon}\text{\ \ by Chebyshev's inequality}\\
&\leq&\frac{\gamma E_{\bm\xi}\left[\sup_{\rho_0(x_1x_2)<\delta}|\tilde G(x_1)-\tilde G(x_2)|\right]}{\epsilon}\\
&\to&0
\end{eqnarray*}
as $\delta\to0$. We have therefore proved \eqref{asymptotic equicontinuity}.
%$$\limsup_{n\to\infty}\sup_{\rho_0(x_1,x_2)<\delta}\widehat{Var}_n(h(x_1;\xi)-h(x_2;\xi))\leq\delta\text{\ \ a.s.}$$
%%where $\delta':=\sup_{\rho_0(x_1,x_2)<\delta}\rho_0(x_1,x_2)\leq\delta$. Since $2\bar\Phi(\cdot)$ is uniformly bounded,
%We then have
%$$\limsup_{n\to\infty}2\bar\Phi\left(\frac{\epsilon}{\sup_{\rho_0(x_1,x_2)<\delta}\left(\widehat{Var}_n(h(x_1;\xi)-h(x_2;\xi))\right)^{1/2}}\right)\leq2\bar\Phi\left(\frac{\epsilon}{\delta}\right)$$
%and hence \eqref{interim new21} equals to 0, concluding that $\tilde G_n(\cdot)$ is uniformly equicontinuous in probability. This,
Together with total boundedness, we have $\tilde G_n(\cdot)\Rightarrow\tilde G(\cdot)$ in $\ell^\infty(\Theta)$ (Section 2.1.2 in \cite{van1996weak}).

Finally, note that $\inf_{x\in\Theta}\sigma_0(x)>0$ by Assumption \ref{nondegeneracy} and Jensen's inequality. Using \eqref{convergence var} and that pointwise division is a continuous map $(\ell^\infty\times\ell^\infty)(\Theta)\to\ell^\infty(\Theta)$, Slutsky's Theorem and the continuous mapping theorem conclude that $G_n(\cdot)\Rightarrow G(\cdot)$ in $\ell^\infty(\Theta)$.
\Halmos
%(CHECK WEAK CONVERGENCE IN GENERAL SPACE)
\endproof

%\begin{corollary}
%Under Assumptions \ref{finite mean}, \ref{nondegeneracy}, \ref{L2} and \ref{complexity}, and selecting $q_n$ as in Theorem \ref{main}, we have, conditional on almost every data realization $(P_n:n\geq1)$,
%$$P(\underline Z_n(x)\leq Z_0(x)\leq\overline Z_n(x)\text{\ for all\ }x\in\Theta)\to1-\alpha$$\label{range cor}
%where
%\begin{align*}
%\underline Z_n(x)&=\min\left\{\sum_{i=1}^nh(x;\xi_i)w_i:-\frac{1}{n}\sum_{i=1}^n\log(nw_i)\leq\frac{q_n}{2n}\ \sum_{i=1}^nw_i=1,\ w_i\geq0\text{\ for\ }i=1,\ldots,n\right\}\\
%\overline Z_n(x)&=\max\left\{\sum_{i=1}^nh(x;\xi_i)w_i:-\frac{1}{n}\sum_{i=1}^n\log(nw_i)\leq\frac{q_n}{2n}\ \sum_{i=1}^nw_i=1,\ w_i\geq0\text{\ for\ }i=1,\ldots,n\right\}
%\end{align*}
%\end{corollary}

\proof{Proof of Theorem \ref{main}.}
By Theorem \ref{EL process thm} and Lemma \ref{sample Gaussian}, and using the fact that $(\cdot)^2$ and $\sup_{x\in\Theta}\cdot$ are continuous maps $\ell^\infty(\Theta)\to\ell^\infty(\Theta)$ and $\ell^\infty(\Theta)\to\mathbb R$ respectively, we have $\sup_{x\in\Theta}\{-2\log R(x;Z_0)\}\Rightarrow\sup_{x\in\Theta}J(x)$ and $\sup_{x\in\Theta}J_n(x)\Rightarrow\sup_{x\in\Theta}J(x)$, where $J_n(\cdot)$ and $J(\cdot)$ are defined in Theorems \ref{main} and \ref{EL process thm}. Moreover, since $\sup_{x\in\Theta}J(x)$ has a continuous distribution function, pointwise convergence of distribution functions to that of $\sup_{x\in\Theta}J(x)$ implies uniform convergence. Hence we have
\begin{equation}
\sup_{q\in\mathbb R^+}\left|P_{\bm\xi}\left(\sup_{x\in\Theta}J_n(x)\leq q\right)-P\left(\sup_{x\in\Theta}J(x)\leq q\right)\right|\to0\label{interim new24}
\end{equation}
and
\begin{equation}
\sup_{q\in\mathbb R^+}\left|P\left(\sup_{x\in\Theta}\{-2\log R(x;Z_0)\}\leq q\right)-P\left(\sup_{x\in\Theta}J(x)\leq q\right)\right|\to0\label{interim new25}
\end{equation}
Selecting $q_n$ such that $P_{\bm\xi}\left(\sup_{x\in\Theta}J_n(x)\leq q_n\right)=1-\alpha$, \eqref{interim new24} and \eqref{interim new25} implies that
$$P\left(\sup_{x\in\Theta}\{-2\log R(x;Z_0)\}\leq q_n\right)\to1-\alpha$$
By applying Proposition \ref{duality} to every point $x\in\Theta$ and with $\chi^2_{1,1-\alpha}$ with $q_n$, we have $-2\log R(x;Z_0)\leq q_n$ if and only if $\underline Z_n^*(x)\leq Z_0(x)\leq\overline Z_n^*(x)$, for each $x\in\Theta$.
%, where
%\begin{align*}
%\underline Z_n(x;q)&=\min\left\{\sum_{i=1}^nh(x;\xi_i)w_i:-\frac{1}{n}\sum_{i=1}^n\log(nw_i)\leq\frac{q}{2n}\ \sum_{i=1}^nw_i=1,\ w_i\geq0\text{\ for\ }i=1,\ldots,n\right\}\\
%\overline Z_n(x;q)&=\max\left\{\sum_{i=1}^nh(x;\xi_i)w_i:-\frac{1}{n}\sum_{i=1}^n\log(nw_i)\leq\frac{q}{2n}\ \sum_{i=1}^nw_i=1,\ w_i\geq0\text{\ for\ }i=1,\ldots,n\right\}
%\end{align*}
Hence
\begin{align*}
P\left(\sup_{x\in\Theta}\{-2\log R(x;Z_0)\}\leq q_n\right)&=P\left(-2\log R(x;Z_0)\leq q_n\text{\ for all\ }x\in\Theta\right)\\%P\left(\{Z(x)\}_{x\in\Theta}\in\left\{\{\tilde Z(x)\}_{x\in\Theta}:\sup_{x\in\Theta}\{-2\log R(\tilde Z(x))\}\leq\eta\right\}\right)=
&=P(\underline Z_n^*(x)\leq Z_0(x)\leq\overline Z_n^*(x)\text{\ for all\ }x\in\Theta)\to1-\alpha
\end{align*}\Halmos
\endproof

\proof{Proof of Theorem \ref{consistency}.}
%Consider
%\begin{equation}
%\max_{\mathbf w\in\hat{\mathcal P}}\sum_{i=1}^nh(x;\xi_i)w_i\label{opt consistency}
%\end{equation}
%where $\hat{\mathcal P}$ is defined as
%$$\hat{\mathcal P}=\left\{\mathbf w=(w_1,\ldots,w_n):-\frac{1}{n}\sum_{i=1}^n\log(nw_i)\leq\eta_n,\ \sum_{i=1}^nw_i=1,\ w_i\geq0\text{\ for\ }i=1,\ldots,n\right\}$$
%where $\eta_n=q_n/(2n)$, and $q_n$ is defined as in Theorem \ref{main}. Let $Z_n(x)$ be the optimal value of \eqref{opt consistency}
To avoid repetition, we focus on $\overline Z_n(x)$. Consider
\begin{equation}
\overline Z_n(x)-Z_0(x)=\max_{\mathbf w\in\mathcal U_n(q_n/(2n))}\sum_{i=1}^n\tilde h(x;\xi_i)w_i\label{opt consistency1}
\end{equation}
where $\tilde h(x;\xi)=h(x;\xi)-Z_0(x)$. With Lagrangian relaxation, the program \eqref{opt consistency1} can be written as
\begin{eqnarray}
&&\min_{\lambda\geq0,\gamma}\max_{\mathbf w\geq\mathbf 0}\sum_{i=1}^n\tilde h(x;\xi_i)w_i-\lambda\left(-\frac{1}{n}\sum_{i=1}^n\log(nw_i)-\frac{q_n}{2n}\right)+\gamma\left(\sum_{i=1}^nw_i-1\right)\notag\\
&=&\min_{\lambda\geq0,\gamma}\sum_{i=1}^n\frac{\lambda}{n}\max_{w_i\geq0}\left\{\frac{\tilde h(x;\xi_i)+\gamma}{\lambda}nw_i+\log(nw_i)-nw_i+1\right\}+\lambda\frac{q_n}{2n}-\gamma\notag\\
&=&\min_{\lambda\geq0,\gamma}-\sum_{i=1}^n\frac{\lambda}{n}\log\left(1-\frac{\tilde h(x;\xi_i)+\gamma}{\lambda}\right)+\lambda\frac{q_n}{2n}-\gamma\label{interim new26}
\end{eqnarray}
where $-0\log(1-t/0):=0$ for $t\leq0$ and $-0\log(1-t/0):=\infty$ for $t>0$, by using the conjugate function of $-\log r+r-1$ as $\sup_{r\geq0}\{tr+\log r-r+1\}=-\log(1-t)$ for $t<1$, and $\infty$ for $t\geq1$ (e.g., \cite{ben2013robust}).

Now, to get an upper bound for \eqref{interim new26}, pick $\gamma=0$, and $\lambda$ as $\lambda_n=\Theta(n^\varepsilon)$ where $1/2<\varepsilon<1$. Then, by using \eqref{interim new1}, we have
$$\max_{1\leq i\leq n}\|\tilde h(\cdot;\xi_i)\|_\Theta\leq\frac{\lambda_n}{2}\text{\ \ ev.}$$
Using the fact that $-\log(1-t)\leq t+2t^2$ for any $|t|\leq1/2$, we have \eqref{interim new26} bounded from above by
\begin{eqnarray}
&&\sum_{i=1}^n\frac{\lambda_n}{n}\left(\frac{\tilde h(x;\xi_i)}{\lambda_n}+2\frac{\tilde h(x;\xi_i)^2}{\lambda_n^2}\right)+\lambda_n\frac{q_n}{2n}\notag\\
&=&\frac{1}{n}\sum_{i=1}^n\tilde h(x;\xi_i)+2\frac{(1/n)\sum_{i=1}^n\tilde h(x;\xi_i)^2}{\lambda_n}+\lambda_n\frac{q_n}{2n}\label{interim new27}
\end{eqnarray}
where $\frac{1}{n}\sum_{i=1}^n\tilde h(x;\xi_i)\to0$ and $\frac{1}{n}\sum_{i=1}^n\tilde h(x;\xi_i)^2$ satisfy
$$\left\|\frac{1}{n}\sum_{i=1}^n\tilde h(\cdot;\xi_i)\right\|_\Theta\to0\text{\ \ a.s.}$$
$$\left\|\frac{1}{n}\sum_{i=1}^n\tilde h(\cdot;\xi_i)^2-\sigma_0^2(\cdot)\right\|_\Theta\to0\text{\ \ a.s.}$$
by the $P_0$-GC property of $\mathcal H_\Theta$ and \eqref{convergence var}. Moreover, by \eqref{interim new24} in the proof of Theorem \ref{main}, we have $P(\sup_{x\in\Theta}J(x)\leq q_n)\to1-\alpha$ a.s.. By the continuity of $\sup_{x\in\Theta}J(x)$, we get $q_n\to q^*$ a.s. where $q^*$ satisfies $P(\sup_{x\in\Theta}J(x)\leq q^*)=1-\alpha$. Hence $q_n/(2n)=\Theta(1/n)$. These imply that \eqref{interim new27} converges to 0 uniformly over $\Theta$.

On the other hand, plugging in $\mathbf w=(1/n)_{1\leq i\leq n}$, $\overline Z_n^*(x)-Z_0(x)$ in \eqref{opt consistency1} is bounded from below by $(1/n)\sum_{i=1}^n\tilde h(x;\xi_i)$, which converges to 0 uniformly over $\Theta$. Combining with above, we get
\begin{equation}
\|\overline Z_n^*(\cdot)-Z_0(\cdot)\|_\Theta\to0\text{\ \ a.s.}\label{interim new28}
\end{equation}

%Finally, consider
%$$Z_0(x_n^*)-Z^*=(Z_0(x_n^*)-Z_n(x_n^*))+(Z_n(x_n^*)-Z_n(x^*))+(Z_n(x^*)-Z_0(x^*))$$
%We have $Z_0(x_n^*)-Z_n(x_n^*)\to0$ and $Z_n(x^*)-Z_0(x^*)\to0$ by \eqref{interim new28}. Moreover, $Z_n(x_n^*)\leq Z_n(x^*)$ by the definition of $x_n^*$. Hence $\limsup_{n\to\infty}Z_0(x_n^*)-Z^*\leq0$ a.s.. Since $Z_0(x_n^*)\geq Z^*$ by the definition of $Z^*$, we have $Z_0(x_n^*)\to Z^*$ a.s..
%
\Halmos
\endproof

\proof{Proof of Theorem \ref{consistency simple}.}
For any fixed $x$, $\underline Z_n(x)\to Z_0(x)$ and $\overline Z_n(x)\to Z_0(x)$ follows as a special case of \eqref{interim new28}.\Halmos
\endproof

\proof{Proof of Theorem \ref{asymptotic equivalence}.}
To avoid redundancy, we focus only on the upper bound $\overline Z_n^*(x)$. Consider $\overline Z_n^*(x)-\hat h(x)$, which can be written as
\begin{equation}
\max\left\{\sum_{i=1}^nw_i\hat h(x;\xi_i):-\frac{1}{n}\sum_{i=1}^n\log(nw_i)\leq\frac{q_n}{2n},\ \sum_{i=1}^nw_i=1,\ w_i\geq0\text{\ for all\ }i=1,\ldots,n\right\}\label{interim asymptotic}
\end{equation}
where $\hat h(x;\xi_i)=h(x;\xi_i)-\hat h_n(x)$. Similar to the proof of Theorem \ref{consistency}, a Lagrangian relaxation of \eqref{interim asymptotic} gives
\begin{eqnarray}
&&\min_{\lambda\geq0,\gamma}\max_{\mathbf w\geq\mathbf 0}\sum_{i=1}^nw_i\hat h(x;\xi_i)-\lambda\left(-\frac{1}{n}\sum_{i=1}^n\log(nw_i)-\frac{q_n}{2n}\right)+\gamma\left(\sum_{i=1}^nw_i-1\right)\notag\\
&\leq&\min_{\lambda>0,\gamma}\sum_{i=1}^n\frac{\lambda}{n}\max_{w_i\geq0}\left\{nw_i\frac{\hat h(x;\xi_i)+\gamma}{\lambda}+\log(nw_i)-nw_i+1\right\}+\frac{\lambda q_n}{2n}-\gamma\notag\\
&\leq&\min_{\lambda>0,\gamma,\frac{\hat h(x;\xi_i)+\gamma}{\lambda}<1\text{\ for all\ }i=1,\ldots,n}-\sum_{i=1}^n\frac{\lambda}{n}\log\left(1-\frac{\hat h(x;\xi_i)+\gamma}{\lambda}\right)+\frac{\lambda q_n}{2n}-\gamma\label{interim asymptotic1}
\end{eqnarray}
by using the fact that $-\log(1-t)=\sup_{r\geq0}\{tr+\log r-r+1\}$ is the conjugate function of $-\log r+r-1$, defined for $t<1$.

Now, given $x$, we choose $\gamma=0$, and $\lambda=\frac{\sqrt n\hat\sigma(x)}{\sqrt{q_n}}$. Similar to the proof of Theorem \ref{consistency}, we have
$$\left\|\frac{1}{n}\sum_{i=1}^n\hat h(\cdot;\xi_i)^2-\sigma_0^2(\cdot)\right\|_\Theta\to0\text{\ \ a.s.}$$
by the $P_0$-GC property of $\mathcal H_\Theta$ and \eqref{convergence var}, and $q_n\to q^*$ a.s. where $q^*$ satisfies $P(\sup_{x\in\Theta}J(x)\leq q^*)=1-\alpha$. Moreover,
$$\hat\sigma^2(x)=\frac{1}{n}\sum_{i=1}^n(h(x;\xi_i)-Z_0(x))^2-(\hat h(x)-Z_0(x))^2$$
By Assumptions \ref{finite mean}, \ref{L2} and \ref{complexity}, $\{(h(x;\cdot)-Z_0(x))^2:x\in\Theta\}$ is $P_0$-Donsker, and we have
$$\sup_{x\in\Theta}|\hat\sigma^2(x)-\sigma_0^2(x)|\to0\text{\ \ as $n\to\infty$ a.s.}$$
Together with the assumption that $h$ is bounded, we have
\begin{equation}
\frac{\hat h(x;\xi_i)+\gamma}{\lambda}=\frac{\hat h(x;\xi_i)\sqrt{q_n}}{\sqrt n\hat\sigma(x)}\to0\label{interim asymptotic2}
\end{equation}
a.s. uniformly for all $i=1,\ldots,n$ as $n\to\infty$.

Therefore, \eqref{interim asymptotic1} is bounded from above ev. by
\begin{eqnarray}
&&-\sum_{i=1}^n\frac{\hat\sigma(x)}{\sqrt{nq_n}}\log\left(1-\frac{\hat h(x;\xi_i)\sqrt q_n}{\sqrt n\hat\sigma(x)}\right)+\frac{\sqrt q_n\hat\sigma(x)}{2\sqrt n}\notag\\
&=&\sum_{i=1}^n\frac{\hat\sigma(x)}{\sqrt{nq_n}}\left(\frac{\tilde h(x;\xi_i)\sqrt q_n}{\sqrt n\hat\sigma(x)}+\frac{1}{2}\left(\frac{\tilde h(x;\xi_i)\sqrt q_n}{\sqrt n\hat\sigma(x)}\right)^2+O\left(\left(\frac{\hat h(x;\xi_i)\sqrt q_n}{\sqrt n\hat\sigma(x)}\right)^3\right)\right)+\frac{\sqrt q_n\hat\sigma(x)}{2\sqrt n}{}\notag\\
&&{}\text{\ \ where $O(\cdot)$ is uniform over $x\in\Theta$}\notag\\
&=&\frac{\sqrt{q_n}\hat\sigma(x)}{2\sqrt n}+O\left(\frac{\hat\mu_3(x)q_n}{n\hat\sigma(x)^2}\right)+\frac{\sqrt q_n\hat\sigma(x)}{2\sqrt n}{}\notag\\
&&{}\text{\ \ since $\sum_{i=1}^n\hat h(x;\xi_i)=0$ and $\frac{1}{n}\sum_{i=1}^n\hat h(x;\xi_i)^2=\hat\sigma^2(x)$,}{}\notag\\
&&{}\text{\ \ where $\hat\mu_3(x)=\frac{1}{n}\sum_{i=1}^n|\hat h(x;\xi_i)|^3$, which is uniformly bounded over $x\in\Theta$ a.s. since $h$ is bounded}\notag\\
&=&\frac{\sqrt{q_n}\hat\sigma(x)}{\sqrt n}+O\left(\frac{1}{n}\right)\label{interim asymptotic3}
\end{eqnarray}

On the other hand, choose
$$w_i=\frac{1}{n}\left(1+\frac{\hat h(x;\xi_i)\sqrt{q_n}}{\sqrt n\hat\sigma(x)}\left(1-\frac{C}{\sqrt n}\right)\right)$$
for some large enough $C>0$. When $n$ is large enough, we have $w_i>0$ a.s. for all $i=1,\ldots,n$ and $x\in\Theta$ since $\frac{\hat h(x;\xi_i)\sqrt{q_n}}{\sqrt n\hat\sigma(x)}<1$ ev. by the same argument in \eqref{interim asymptotic2}. Note that $\sum_{i=1}^nw_i=1$ since $\sum_{i=1}^n\hat h(x;\xi_i)=0$ by definition. Moreover,
\begin{eqnarray*}
&&-\frac{1}{n}\sum_{i=1}^n\log(nw_i)\\
&=&-\frac{1}{n}\sum_{i=1}^n\log\left(1+\frac{\hat h(x;\xi_i)\sqrt{q_n}}{\sqrt n\hat\sigma(x)}\left(1-\frac{C}{\sqrt n}\right)\right)\\
&=&-\frac{1}{n}\sum_{i=1}^n\frac{\hat h(x;\xi_i)\sqrt{q_n}}{\sqrt n\hat\sigma(x)}\left(1-\frac{C}{\sqrt n}\right)+\frac{1}{n}\frac{1}{2}\sum_{i=1}^n\frac{\hat h(x;\xi_i)^2q_n}{n\hat\sigma(x)^2}\left(1-\frac{C}{\sqrt n}\right)^2+O\left(\frac{\hat\mu_3(x)q_n^{3/2}}{n^{3/2}\hat\sigma(x)^3}\left(1-\frac{C}{\sqrt n}\right)^3\right){}\\
&&{}\text{\ \ where $O(\cdot)$ is uniform over $x\in\Theta$}\\
&=&\frac{q_n}{2n}\left(1-\frac{2C}{\sqrt n}\right)+O\left(\frac{1}{n^2}\right)+O\left(\frac{\hat\mu_3(x)q_n^{3/2}}{n^{3/2}\hat\sigma(x)^3}\right){}\\
&&{}\text{\ \ where the last $O(\cdot)$ has leading term that is independent of $C$}\\
&\leq&\frac{q_n}{2n}
\end{eqnarray*}
when $n$ is large enough, by choosing a large enough $C$. Therefore, the chosen $w_i$'s form a feasible solution in $\mathcal U_n(q_n/(2n))$. We have
\begin{align}
\sum_{i=1}^nw_i\hat h(x;\xi_i)&=\sum_{i=1}^n\hat h(x;\xi_i)\frac{1}{n}\left(1+\frac{\hat h(x;\xi_i)\sqrt{q_n}}{\sqrt n\hat\sigma(x)}\left(1-\frac{C}{\sqrt n}\right)\right)\notag\\
&=\sqrt{q_n}\frac{\hat\sigma(x)}{\sqrt n}\left(1-\frac{C}{\sqrt n}\right)\notag\\
&=\sqrt{q_n}\frac{\hat\sigma(x)}{\sqrt n}+O\left(\frac{1}{n}\right)\label{interim asymptotic4}
\end{align}

Combining the bound for the dual and the primal bounds \eqref{interim asymptotic3} and \eqref{interim asymptotic4}, we conclude that $\overline Z_n^*(x)=\sqrt{q_n}\frac{\hat\sigma(x)}{\sqrt n}+O\left(\frac{1}{n}\right)$ uniformly over $x\in\Theta$. The proof for $\underline Z_n^*(x)$ follows by merely replacing $h$ with $-h$. This concludes the theorem.\Halmos
\endproof

\section{Review of Empirical Processes}\label{sec:EP}
%Consider the random object $\xi\in\Xi$ under probability measure $P$, and expectation is denoted by $E$. Consider the function class $\mathcal F$ that maps from $\mathcal X$ to $\mathbb R$. Given i.i.d. $X_1,\ldots,X_n$, define the empirical measure $\mathbb P_n:\mathcal F\to\mathbb R$ as
We review some terminologies and results in the empirical process theory that are related to our developments. %To begin with, let us denote the class of functions mapping from $\Xi$ to $\mathbb R$, indexed by $\Theta$,
%$$\mathcal F_\Theta:=\{f(x;\cdot):\Xi\to\mathbb R|x\in\Theta\}$$
Given a class of functions $\mathcal F=\{f:\Xi\to\mathbb R\}$, we define the empirical measure $\mathbb P_n$, generated from i.i.d. $\xi_1,\ldots,\xi_n$ each under $P$, as a map from $\mathcal F$ to $\mathbb R$ such that
$$\mathbb P_n(f)=\frac{1}{n}\sum_{i=1}^nf(\xi_i)$$
%The empirical measure indexed by $\Theta$ refers to the collection of empirical measures
%$$\{\mathbb P_n(x):x\in\Theta\}$$
We also define $P(f)=\int f(\xi)dP(\xi)=E_P[f(\xi)]$ where $E_P[\cdot]$ is the expectation under $P$. The empirical process indexed by $f\in\mathcal F$ is defined as
$$\sqrt n(\mathbb P_n-P)$$
%The empirical process indexed by $\Theta$ refers to the collection of empirical processes
%$$\{\sqrt n(\mathbb P_n(x)-P(x)):x\in\Theta\}$$
For any functions $y:\mathcal F\to\mathbb R$, we define $\|y\|_{\mathcal F}=\sup_{f\in\mathcal F}|y(f)|$. We also define the envelope of $\mathcal F$ as a function that maps from $\Xi$ to $\mathbb R$ given by
$$\sup_{f\in\mathcal F}|f(\xi)|$$

\begin{definition}
We call $\mathcal F$ a $P$-Glivenko-Cantelli (GC) class if the empirical measure under $P$ satisfies
$$\|\mathbb P_n-P\|_{\mathcal F}:=\sup_{f\in\mathcal F}|\mathbb P_n(f)-P(f)|\stackrel{a.s.}{\to}0\text{\ \ as\ }n\to\infty$$
\end{definition}
\begin{definition}
We call $\mathcal F$ a $P$-Donsker class if the empirical process under $P$ satisfies
\begin{equation}
\sqrt n(\mathbb P_n-P)\Rightarrow\mathbb G\text{\ \ in\ }\ell^\infty(\mathcal F)\label{Donsker thm}
\end{equation}
where $\mathbb G$ is a Gaussian process indexed by $\mathcal F$, centered, with covariance function
$$Cov(\mathbb G(f_1),\mathbb G(f_2))=Cov_P(f_1(\xi),f_2(\xi))=P(f_1f_2)-P(f_1)P(f_2)$$
where $Cov_P(\cdot,\cdot)$ denotes the covariance under $P$, and
$$\ell^\infty(\mathcal F)=\left\{y:\mathcal F\to\mathbb R\Bigg|\|y\|_{\mathcal F}<\infty\right\}$$
Moreover, the process $\mathbb G$ has uniformly continuous sample paths with respect to the canonical semi-metric $\rho_P(f_1,f_2)=Var_P(f_1(\xi)-f_2(\xi))$, where $Var_P(\cdot)$ denotes the variance under $P$.\label{def:Donsker}
\end{definition}

We have ignored the measurability issues, in particular the use of outer and inner probability measures, in the definitions (see \cite{van1996weak}). %, in that the space of all continuous functions $\ell^\infty(\mathcal F)\to\mathbb R$ required to define the weak convergence can be non-separable.
%For detailed discussion, we refer the reader to \cite{van1996weak}.

%As a sanity check, note that $\{\chi^2(x)\}\in\ell^\infty(\Theta)$ since $\Theta$ is separable (REF).
%Consider a cost function $h(X;\theta)$ with decision $\theta\in\Theta$, where $X\in\mathcal X$ is a random object under $P$. The goal is to find the optimal solution $\theta^*$ of the stochastic optimization $\min_{\theta\in\Theta}Z(\theta)=E[h(X;\theta)]$. Suppose we have i.i.d. data $X_1,\ldots,X_n$. We posit the distributionally robust optimization (DRO)

%We denote $\hat Z_n(\theta)=\max_{\mathbf w\in\mathcal W}\sum_{i=1}^nh(X_i;\theta)w_i$, and $\hat Z_n^*$ the optimal value of \eqref{opt}. Thus $\hat Z_n^*=\min_\theta\hat Z(\theta)$. Let $\hat\theta_n^*$ be the optimal solution to \eqref{opt}.
%
%We sometimes suppress $\theta$ by writing $h(X)=h(X;\theta)$ and $Z=Z(\theta)$.
%
%Note that $\mathbb G$ is totally bounded.

\begin{theorem}[Preservation of GC classes; \cite{van2000preservation}, Theorem 3]
Suppose that $\mathcal F_1,\ldots,\mathcal F_k$ are $P$-GC classes of functions, and that $\varphi:\mathbb R^k\to\mathbb R$ is continuous. Then $\mathcal H=\varphi(\mathcal F_1,\ldots,\mathcal F_k)$ is $P$-GC provided that it has an integrable envelope function.\label{preservation}
\end{theorem}

\section{Other Useful Theorems}\label{sec:thms}
\begin{lemma}[\cite{owen2001empirical}, Lemma 11.2]
Let $Y_i$ be i.i.d. random variables on $\mathbb R$ with $EY_i^2<\infty$. Then $\max_{1\leq i\leq n}|Y_i|=o(n^{1/2})$ a.s..\label{lemma max}
\end{lemma}

\begin{theorem}[Continuous Mapping Theorem;\cite{van1996weak},Theorem 1.3.6]
Let $g:\mathbb D\to\mathbb E$ be continuous at every point $\mathbb D_0\subset\mathbb D$. If $X_n\Rightarrow X$ and $X$ takes its values in $\mathbb D_0$, then $g(X_n)\Rightarrow g(X)$.\label{continuous mapping}
\end{theorem}

\begin{theorem}[Slutsky's Theorem; \cite{van1996weak}, Example 1.4.7]
If $X_n\Rightarrow X$ and $Y_n\Rightarrow c$ where $X$ is separable and $c$ is a constant, then $(X_n,Y_n)\Rightarrow(X,c)$ under the product topology.\label{Slutsky}
\end{theorem}

\end{APPENDICES}

\end{document}